\documentclass[10pt]{smfart}
\usepackage{a4wide}
\usepackage{amsmath,amssymb,smfthm}
\usepackage[francais]{babel}
\usepackage[latin1]{inputenc}
\usepackage[all]{xy}
\usepackage{graphicx}
\usepackage{epic,eepic}
\usepackage{epsfig}
\usepackage{color}

\def\R{\mathbb{R}}

\def\P{\mathbb{P}}
\def\C{\mathbb{C}}

\def\Cp{\mathbb{C}_p}
\def\Fp{\mathbb{F}_p}

\def\Fqb{\overline{\mathbb{F}}_q}
\def\L{\mathcal{L}}

\def\O{\mathcal{O}}

\def\*{^\times }
\def\dpt{\displaystyle}

\def\l{\lambda}

\def\a{\alpha}

\def\s{\sigma}
\def\ph{\varphi}

\def\e{\epsilon}
\def\ssi{\Leftrightarrow}

\def\impl{\Rightarrow}

\def\drt{\rightarrow}
\def\ldrt{\longrightarrow}
\def\Q{\mathbb{Q}}

\def\Qp{\mathbb{Q}_p}
\def\Qpb{\overline{\mathbb{Q}}_p}
\def\Zp{\mathbb{Z}_p}
\def\Z{\mathbb{Z}}
\def\N{\mathbb{N}}

\def\Hom{\text{Hom}}

\def\Gal{\text{Gal}}
\def\={\! = \!}

\def\spec{\text{Spec}}
\def\spf{\text{Spf}}

\def\limp{\underset{\longleftarrow}{\text{ lim }}\;}

\def\iso{\xrightarrow{\;\sim\;}}

\def\End{\text{End}}

\def\GL{\hbox{GL}}

\def\xrig{\xrightarrow}

\def\M{\mathcal{M}}

\def\X{\mathfrak{X}}
\def\GG{\Gamma}

\def\bc{\backslash}
\def\AA{\mathcal{A}}

\def\<{<\hspace{-1mm}}
\def\>{\hspace{-1mm}>}

\def\Lie{\text{Lie}}

\def\dem{{\it Démonstration. }}

\def\Fil{\text{Fil}}

\def\unp{[ \frac{1}{p}]}
\def\LT{\mathcal{L}\mathcal{T}}
\def\D{{\mathcal{D}r}}
\def\Hb{\mathbb{H}}
\def\Gb{\mathbb{G}}
\def\DH{\mathbb{D}  (\mathbb{H}  )}
\def\DG{\mathbb{D} (\mathbb{G}  )}

\author{Laurent Fargues}
\address{CNRS-IHES-université Paris-Sud Orsay}
\email{laurent.fargues@math.u-psud.fr}

\begin{document}

\newtheorem{Fait}{Fait}

\title{L'isomorphisme entre les tours de Lubin-Tate et de Drinfeld au niveau des points}
\maketitle

\begin{abstract}
Dans cet article on construit l'isomorphisme entre les tours de
Lubin-Tate et de Drinfeld au niveau des points de ces espaces. Les
points considérés sont ceux intervenant dans la théorie des espaces
analytiques de Berkovich.
\end{abstract}

\begin{altabstract}
In this article we construct the isomorphism between Lubin-Tate and
Drinfeld towers at the level of points. The points we consider are the
one of the theory of analytic spaces in the sens of Berkovich. 
\end{altabstract}

\section*{Introduction}

Le but de cet article est de démontrer l'existence de l'isomorphisme entre les tours
de Lubin-Tate et de Drinfeld au niveau des points c'est à dire de
décrire la bijection correspondante. Nous nous inspirons bien sûr de \cite{Faltings8}
(cependant l'auteur ne garantit pas que la bijection décrite coïncide
avec celle de \cite{Faltings8} qu'il n'a pu complètement comprendre). 

Cet isomorphisme n'est pas algébrique : un $\Qpb$-point peut s'envoyer
sur un $\Cp$-point ne provenant pas d'un $\Qpb$-point. C'est pourquoi
nous devons travailler avec des points à valeurs dans de ``gros
corps'' du type $\Cp$. 
Par points nous entendons donc ceux intervenant par exemple dans la
théorie des espaces analytiques de Berkovich : ce sont ceux à valeurs
dans un corps valué complet pour une valuation de rang $1$ extension
de $\Qp$.

La simplification par rapport à la construction de l'isomorphisme dans
le cas général provient de ce que l'on n'a pas à introduire
d'éclatements admissibles ou de modèles entiers particuliers des
espaces de Lubin-Tate et de Drinfeld : si $K|\Qp$ est valué complet
pour une valuation de rang $1$ les points à valeurs dans $K$ de la
tour de Lubin-Tate et de Drinfeld sont définis indépendamment du
modèle entier. En particulier cet article est indépendant
de la construction du schéma formel de  \cite{Cellulaire}.

L'utilisation de corps valués pour des valuations non-discrètes à
corps résiduel non forcément parfait introduit des problèmes de
théorie de Hodge $p$-adique non disponibles dans la littérature. Les
résultats associés peuvent se déduire de l'approche utilisée dans
\cite{Faltings1}. Nous ne les démontrons pas dans cet article mais
nous y consacrons un appendice dans lequel nous y exposons les
résultats auxquels on peut parvenir en utilisant les méthodes de
\cite{Faltings1}. L'auteur y consacrera un article plus général
(\cite{Periodes}). 

L'un des autres problèmes auquel on est confronté est le fait que l'on
travaille avec des $\O$-modules formels, le cas usuel correspondant à $\O=\Zp$.
 Pour y remédier on applique
la théorie des $\O$-extensions vectorielles universelles 
développée dans l'appendice B de \cite{Cellulaire}. Très peu est
nécessaire : seules les sections un à quatre de ce cet appendice sont
utilisées, nous n'utiliserons pas le fait que le cristal algèbre de
Lie de la $\O$-extension vectorielle universelle s'étend aux
$\O$-puissances $\pi$-divisées. qui est la partie délicate de
l'appendice B de \cite{Cellulaire}. L'utilisation de cette théorie
relative à $\O$ rend la rédaction des problèmes reliés purement
formels. Seul le cas de l'utilisation des théorèmes de comparaison
cristallins n'est pas rédigé dans ce cadre relatif : l'auteur n'a pas
eu le courage de développer la théorie des anneaux de Fontaine obtenus
en remplaçant les vecteurs de Witt par les vecteurs de Witt
ramifiés...

Dans cet article nous donnons deux descriptions différentes de
l'isomorphisme. La première dans le chapitre 9 est la plus simple.
La seconde dans le chapitre 12 est plus alambiquée mais s'adapte mieux
au cas relatif (on renvoie à l'introduction du chapitre \ref{ptovue}
pour plus de détails).
\\

Avant de décrire succinctement la bijection rappelons quelques faits
sur les espaces de Lubin-Tate et de Drinfeld
associés à $F|\Qp$. 
 Plaçons nous au
niveau des points de ces espaces. 
Soit $D$ l'algèbre à division
d'invariant $\frac{1}{n}$ sur $F$.
La tour de Lubin-Tate $(\LT_K)_{K\subset GL_n (\O_F)}$  est une tour
indexée par les sous-groupes ouverts $K$ de $\GL_n (\O_F)$ 
munie d'une action ``horizontale'' de $D^\times$ (sur chaque élément
de la tour) 
 et verticale (par
correspondances de Hecke) de $\GL_n (F)$. Soit $\LT_\infty
=\underset{K}{\limp} \LT_K$ (cela a bien un sens au niveau des points)
qui est muni d'une action de $\GL_n(F)\times D^\times$. Elle forme un
pro-revêtement de groupe $\GL_n (\O_F)$
$$
\xymatrix@C=5mm{
\ar@/_1.2pc/@{-}[d]_{K}
\LT_\infty \ar[d]  \ar@/^2pc/@{-}[dd]^{GL_n (\O_F)}
\\
\LT_K \ar[d] \ar@(l,ld)_{D^\times}
\\
\LT_{GL_n (\O_F)} \ar@{=}[r] &  \coprod_\Z {\mathring{\mathbb{B}}^{n-1}}
}
$$
où $\LT_{GL_n (\O_F)}$ est l'espace de Lubin-Tate sans niveau (une
union disjointes de
boules $p$-adiques ouvertes de dimension $n-1$ car on travaille en
fait avec des espaces de Rapoport-Zink et non les espaces de
Lubin-Tate classiques) qui classifie des
déformations de groupes formels de dimension $1$ et $\LT_K$ est
obtenu en trivialisant partiellement le module de Tate de la
déformation universelle au dessus de $\LT_{GL_n (\O_F)}$ (on force la
monodromie sur le module de Tate à vivre dans $K$). 

On a donc $\LT_\infty/\GL_n (\O_F) = \coprod_\Z\mathring{\mathbb{B}}^{n-1}$. Il
y a de plus une application des périodes de Hodge-De-Rham $\LT_K \ldrt
\mathbb{P}^{n-1}$ donnée par la filtration de Hodge du module filtré
définissant la déformation et 
dont les fibres sont les orbites de Hecke. Et donc
$\LT_\infty/\GL_n (F)= \mathbb{P}^{n-1}$. 
$$
\xymatrix@R=8mm{
\LT_\infty \ar[d]  \ar@/^2pc/@{-}[d]^{GL_n (\O_F)}
\ar@/_2pc/@{-}[dd]_{GL_n (F)} \\
  {\coprod_\Z \mathring{\mathbb{B}}^{n-1}} \ar[d]^{\text{Périodes}} \\
\mathbb{P}^{n-1} \ar@(l,ld)_{D^\times} 
}
$$
Il y a une description du même type pour la tour de Drinfeld. 
L'espace de Drinfeld en niveau infini $\D_\infty$ est muni d'une
action de $\GL_n (F)\times D^\times$, l'action de $\GL_n (F)$ est
``horizontale'' sur chaque élément de la tour.
Il y  a un
diagramme pour $K\subset \O_D^\times$ un sous-groupe ouvert 
$$
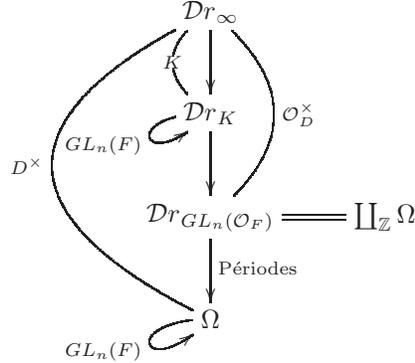
\xymatrix{
\ar@/_1.2pc/@{-}[d]|{K} \ar@/_5pc/@{-}[ddd]_{D^\times}
\D_\infty \ar[d]  \ar@/^2pc/@{-}[dd]^{\O_D^\times}
\\
\D_K \ar[d] \ar@(l,ld)_{GL_n (F)}
\\
\D_{GL_n (\O_F)} \ar@{=}[r]\ar[d]^{\text{Périodes}}  &  \coprod_\Z \Omega \\
\Omega\ar@(l,ld)_{GL_n (F)}
}
$$
L'application des périodes de Hodge-De-Rham n'est rien d'autre que la
projection de $\coprod_\Z \Omega$ sur $\Omega$, c'est à dire le
quotient par $D^\times/\O_D^\times = \Pi^\Z$. 
\\

Les deux applications inverses l'une de l'autre entre $\LT_\infty$ et $\D_\infty$ sont construites en deux
étapes. 
$$
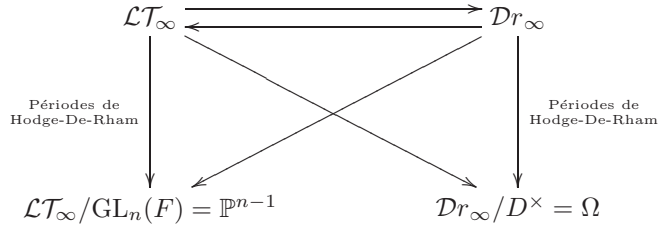
\xymatrix@R=20mm{
\LT_\infty \ar[d]_{\text{Périodes de} \atop{\text{Hodge-De-Rham}}}\ar@<.8ex>[rr] \ar[rrd]  && \D_\infty \ar[d]^{\text{Périodes de} \atop{\text{Hodge-De-Rham}}} \ar[lld] \ar@<.8ex>[ll] \\
\LT_\infty/\GL_n (F) = \mathbb{P}^{n-1} & & \D_\infty/D^\times =\Omega
}
$$
On commence par construire une application de périodes de Hodge-Tate d'un des
espaces en niveau infini vers l'espace des périodes de Hodge-De-Rham
de l'autre espace (les flèches diagonales dans la figure précédente). 
Les périodes de Hodge-Tate sont données par la
rigidification du module de Tate (le fait qu'on ai fixé une base de ce
module au sommet de la tour) et la suite de Hodge-Tate,
 contrairement aux périodes de
Hodge-De-Rham données elles par la rigidification de la cohomologie
cristalline (l'isocristal du groupe $p$-divisible) et la filtration de
Hodge. Puis on relève la
flèche en niveau infini en utilisant la théorie de Messing qui permet
de construire des éléments dans le module de Tate en les construisant
modulo $p$. 

Nous n'expliquerons pas d'avantage ce dernier point dans
l'introduction, mais afin d'en donner un avant-goût  citons le
point clef suivant : soit $K|\Qp$ une extension valuée complète et $H$ un groupe
$p$-divisible sur $\O_K$ . Notons $\bar{\eta}=\spec
(\overline{K})$. Alors le module de Tate de $H$ admet deux définitions
: l'une en fibre générique, la définition usuelle
$$
T_p (H) =\underset{n}{\limp} H[p^n] (\overline{K}) = \Hom (\Qp/\Zp, H_{\bar{\eta}})
$$
l'autre modulo $p$ 
$$
T_p (H) = \{\; f \in \Hom (\Qp/\Zp, H\otimes \O_{\overline{K}}/p
\O_{\overline{K}})\;|\;
f_* (1) \in \Fil\, \Lie (E(H))\otimes \O_{\widehat{\overline{K}}}\; \}
$$
où $\Lie\, E(H)$ désigne l'algèbre de Lie de l'extension vectorielle
universelle de $H$ muni de sa filtration de Hodge 
et 
$f_* : \O_{\widehat{\overline{K}}} \ldrt \Lie \, E(H)\otimes \O_{\widehat{\overline{K}}}$ est le morphisme
induit sur les cristaux de Messing évalués sur l'épaississement à
puissances divisées $\O_{\widehat{\overline{K}}} \twoheadrightarrow \O_{\overline{K}}/p\O_{\overline{K}}$.
La démonstration est un jeu basé sur ces deux aspects : fibre
générique et modulo $p$. 
 Les deux groupes $p$-divisibles sur $\overline{\mathbb{F}}_p$ que l'on
déforme pour définir les espaces de Lubin-Tate et de Drinfeld sont
reliés modulo $p$, l'un est isogéne à un produit d'un nombre
fini de copies de l'autre,  et l'on peut ainsi ``transférer'' des éléments du
module de Tate sur une des tours vers l'autre.
\\

Une autre interprétation de cet isomorphisme est donnée dans la
section \ref{kfyygpzr5rrppn}. On peut associer à un point de l'une des
deux tours en niveau infini une matrice de périodes dans $M_n
(B^+_{cris})$. En effet, pour la tour de Lubin-Tate
le théorème de comparaison entre cohomologie
cristalline et cohomologie étale $p$-adique pour les groupes
$p$-divisibles fournit une matrice de périodes $X_{cris}$ 
 à coefficients dans
$B^+_{cris}$ une fois que l'on a fixé une base du module de Tate (la
rigidification en niveau infini) et une base du module de Dieudonné
(la rigidification définissant la déformation)
$$
(B^+_{cris})^n\iso V_p\otimes B^+_{cris} \iso \mathbb{D}^{\LT}\otimes
{B^+_{cris}} \xleftarrow{\; \sim\;} (B^+_{cris})^n
$$
Pour la tour de Drinfeld il y a un isomorphisme de $D\otimes_{\Qp} B^+_{cris}$-modules
en niveau infini
$$
D\otimes_{\Qp} B^+_{cris} \iso V_p \otimes B^+_{cris} \iso \mathbb{D}^\D\otimes
B^+_{cris} \xleftarrow{\; \sim\;} D\otimes_{\Qp} {B^+_{cris}}
$$
 utilisant la $D$-équivariance de l'isomorphisme des périodes
cristalline on en déduit une matrice $Y_{cris}\in M_n (B^+_{cris})$ comme précédemment.

L'isomorphisme entre
les deux tours consiste alors à prendre la transposé de ces matrices
de périodes
$$
Y_{cris} = \,^t X_{cris}
$$ 
Si $X_{cris}$ désigne la matrice de périodes sur
$\LT_\infty$ sa réduction modulo l'idéal d'augmentation
$B^+_{cris}\ldrt \O_{\widehat{\overline{K}}}$, $X$, est telle que les
colonnes de $X$ engendrent les périodes de Hodge-De-Rham et ses lignes
les périodes de Hodge-Tate !
$$
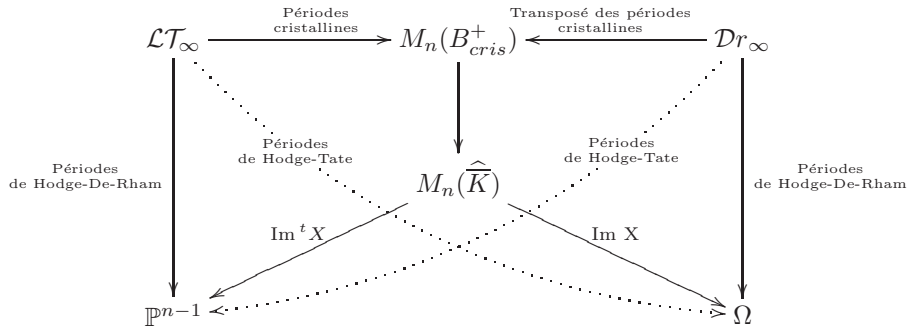
\xymatrix@C=24mm@R=12mm{
\LT_\infty \ar[dd]_{\text{Périodes}\atop\text{de Hodge-De-Rham}}
\ar[r]^{\text{Périodes}\atop \text{cristallines}} 
\ar@/_2pc/@{..>}[rrdd]|(.25){\text{Périodes}\atop\text{de Hodge-Tate}}
& M_n
(B^+_{cris}) \ar[d] & \ar[l]_{\text{Transposé des périodes}\atop \text{cristallines}} \D_\infty
\ar[dd]^{\text{Périodes}\atop\text{de Hodge-De-Rham}}
\ar@/^2pc/@{..>}[lldd]|(.25){\text{Périodes}\atop\text{de Hodge-Tate}} \\
 &  M_n ({\widehat{\overline{K}}}) \ar[ld]_{\text{Im}\,^t
   X} \ar[rd]^{\text{Im X}} \\
\mathbb{P}^{n-1} && \Omega 
}
$$
(les deux flêches $M_n (\widehat{\overline{K}})\ldrt \mathbb{P}^{n-1}$
et $M_n (\widehat{\overline{K}})\ldrt \Omega$ ne sont définies que sur
le sous-ensemble des matrices de rang $n-1$). 
Dans le diagramme précédent 
l'image des deux flèches ``Périodes cristallines'' sont les matrices
$X_{cris}\in M_n (B^+_{cris})$ vérifiant une équation fonctionnelle
``à la Fontaine'' 
du type $$ X_{cris}^{\ph} \Psi = p X_{cris}$$ où $\Psi \in \GL_n (
W(\overline{\mathbb{F}}_p)\unp)$ est la matrice d'un Frobenius sur un
certain isocristal, $X_{cris}^\ph$ est obtenue en appliquant le
Frobenius cristallin aux termes de la matrice 
 et la réduite via $B^+_{cris}\ldrt
\widehat{\overline{K}}$, 
$X\in M_n (\widehat{\overline{K}})$ vérifie 
$$
\text{rang} (X) =n-1
$$
\\

{\it 
Prérequis : La chapitre 3 fournit les rappels nécessaires concernant
les espaces de Lubin-Tate et de Drinfeld. Le seul prérequis
non rappelé est la théorie de la déformation de Grothendieck-Messing
(\cite{Messing1}).  
}
\\

{\it Avertissements : Dans tout l'article on supposera $p\neq 2$. Il
  est conseillé au lecteur de supposer en première lecture que
  $\O=\Zp$.
}
\\

{\it Remerciement : L'auteur remercie Jean François Dat qui,
à partir d'une version préliminaire de cet article, a donné un
  exposé sur le sujet au groupe de travail de l'IHES au printemps 2004
  et ainsi suggéré des améliorations.
}

\section{Décomposition de Hodge-Tate des groupes $p$-divisibles dans le cas infiniment ramifié}\label{HoTa}

Soit $K|\Qp$ un corps valué complet pour une valuation à valeurs dans
$\R$ étendant celle de $\Qp$.
\\

Le théorème suivant se déduit du théorème \ref{bouqibouc} de
l'appendice, lui même démontré dans \cite{Periodes}. 

\begin{theo}[\cite{Periodes}]\label{decHTO}
Soit $G$ un groupe $p$-divisible sur $\O_K$. 
Il y a une décomposition de Hodge-Tate : une suite de $\Gal (\overline{K}|K)$-modules continus 
$$
0\ldrt \omega_G^* \otimes  \mathcal{O}_{\widehat{\overline{K}}} (1) \ldrt T_p (G)\otimes_{\Zp}   \mathcal{O}_{\widehat{\overline{K}}} \ldrt \omega_{G^D} \otimes  \mathcal{O}_{\widehat{\overline{K}}} \ldrt 0
$$
dont la cohomologie comme suite de
$\mathcal{O}_{\widehat{\overline{K}}}$-modules  est annulée par
$p^{\frac{1}{p-1}}$. Si $\O_{K_0} \subset \O_K$ est un anneau de Cohen
 il existe des morphismes $\O_{K_0}$-linéaires Galois équivariantes des deux cotés de la suite exacte tels que composés avec les applications de droite et de gauche on obtienne la multiplication par un élément de valuation $p$-adique $\frac{1}{p-1}$.  
\end{theo}

Bien s\^ur la suite est définie sans recours à la théorie de Hodge
p-adique, celle-ci ne sert que pour démontrer que la cohomologie
de la suite est annulée par tout élément de valuation supérieure ou égale à $\frac{1}{p-1}$. Rappelons en effet que 
$$
\a_G: T_p (G) \ldrt \omega_{G^D}\otimes \O_{\widehat{\overline{K}}}
$$
est défini de la façon suivante : pour $x\in T_p (G) = \Hom (\Qp/\Zp,
G_{/\O_{\overline{K}}})$ celui ci donne par dualité de Cartier un morphisme 
$$
x^D : G_{/\O_{\overline{K}}}^D\ldrt \mu_{p^\infty}
$$
sur $\spec (\O_{\overline{K}})$, donc sur $\spf (\O_{\widehat{\overline{K}}})$.
 Il induit un morphisme 
$$
(x^D)^*:  \O_{\widehat{\overline{K}}} \frac{dT}{T}\ldrt
 \omega_{G^D}\otimes  \O_{\widehat{\overline{K}}}  
$$
et alors
$$
\a_G ( x) = (x^D)^* \frac{dT}{T}
$$
L'application $ \omega_G^* \otimes
\mathcal{O}_{\widehat{\overline{K}}} (1) \ldrt T_p (G)\otimes_{\Zp}
\mathcal{O}_{\widehat{\overline{K}}}$ est définie par dualisation de
$\a_{G^D}$, c'est $\,^t\a_{G^D}(1)$ après identification de $T_p (G^D)$
avec $T_p (G)^*(1)$.

\begin{rema}
Lorsque $K|K_0$ est de degré fini il existe de telles quasi-sections 
$\O_{K}$-linéaires telles que les composées donnent la multiplication
par un élément de valuation $\frac{1}{p-1} + v_p
(\mathcal{D}_{K/K_0})$ où $\mathcal{D}$ désigne la
différente. cf. \cite{Fontaine1}. Cela peut s'obtenir avec les
méthodes précédentes en calculant l'annulateur de $t$ dans un gradué
de $A_{cris}\otimes_{\O_{K_0}}\O_{K}$. Dans notre cas infiniment
ramifié cette suite exacte n'est en général pas scindée après
inversion de $p$, il n'y a pas de décomposition de Hodge-Tate;
l'obstruction provient de la Connexion de Gauss Manin et du fait que
$\Omega^1_{K/K_0} \neq 0$ alors que ce module est nul  dans l'article
\cite{Fontaine1}. Un tel théorème ne peut s'obtenir par les méthodes
de \cite{Fontaine1} qui sont purement ``de Rham'' : il faut utiliser
les propriétés cristallines de certains objets associés aux groupe $p$-divisibles.
\end{rema}

\begin{rema} \label{HTet}
Si pour une extension $K'$ de $K$ contenue dans $\overline{K}$ 
 l'action de $\Gal( \overline{K} | K)$ sur $T_p (G)$ se factorise via
$\Gal (K' |K)$  
il y a alors une décomposition de Hodge-Tate sur $K'$ c'est à dire une suite 
$$
0\ldrt \omega_G^* \otimes  \mathcal{O}_{\widehat{K'}} (1) \ldrt T_p (G)\otimes_{\Zp}  
 \mathcal{O}_{\widehat{K'}} \ldrt \omega_{G^D} \otimes  \mathcal{O}_{\widehat{K'}} \ldrt 0
$$
dont la cohomologie est annulée par tout élément de valuation $\geq \frac{1}{p-1}$. Cela résulte de ce que la suite peut être définie sans recours à la théorie de Fontaine comme ci-dessus
 et que $\O_{\widehat{\overline{K}}}$ est fidèlement plat sur $\mathcal{O}_{\widehat{K'}}$. 
\end{rema}

Soit $\Lie \, E(G)$ l'algèbre de Lie de l'extension vectorielle
universelle de $G$. Elle possède une filtration donnée par la partie
vectorielle 
$$
0 \ldrt \omega_{G^D} \ldrt \Lie \, E(G) \ldrt \omega_G^* \ldrt 0
$$
Alors,
l'application de Hodge-Tate $\a_G$ se décrit également de la
façon suivante 
$$
\xymatrix@R=2mm@C=5mm{
\Hom (\Qp/\Zp, G) \ar[r] & \Hom_{\text{\it Filtré}} \left ( \Lie \, E(\Qp/\Zp),
\Lie \, E(G)\right ) \ar@{=}[r] & \Hom_{\O_K} (\O_K,\omega_{G^D}) \ar[r] &
\omega_{G^D} \\
&& f \ar@{|->}[r] & f(1) 
}
$$
puisque $\Lie \, E(\Qp/\Zp) =\Fil\, \Lie \, E(\Qp/\Zp) = \O_K$.

\subsection{Décomposition de Hodge-Tate d'un $\O$-module
  $\pi$-divisible}

On reprend les notations de l'appendice B de \cite{Cellulaire}. On y
note $F$ une extension de degré fini de $\Qp$ et $\O=\O_F$ son anneau
des entiers. 

Soit $K$ comme dans les sections précédentes et supposons que $K|F$.
Soit $G$ un $\O$-module $\pi$-divisible sur $\O_K$,
c'est à dire un groupe $p$-divisible muni d'une action de $\O$
induisant l'action naturelle sur son algèbre de Lie. 
 Soit la suite
exacte de Hodge-Tate définie dans la section précédente 
$$
0\ldrt \omega_G^*\otimes\widehat{\overline{K}} \ldrt V_p
(G)\otimes_{\Qp} \widehat{\overline{K}} \ldrt \omega_{G^D}
\otimes\widehat{\overline{K}} \ldrt 0
$$
Il s'agit d'une suite de $\O\otimes_{\Zp} \O_K$-modules. 
\\

Rappelons (\cite{Cellulaire}) que l'on pose 
$$
\widetilde{\omega}_{G^D}= \omega_{G^D}/I.\Lie E(G) \text{ où } I=\ker
(\O\otimes_{\Zp} \O_K \twoheadrightarrow \O_K )
$$
qui est la partie vectorielle de la $\O$-extension vectorielle
universelle de $G$. Rappelons en effet que si 
$$
0\ldrt V(G) \ldrt E(G) \ldrt G \ldrt 0
$$
est l'extension vectorielle universelle de $G$ alors la $\O$-extension
vectorielle universelle est le poussé en avant par l'application
$V(G)\twoheadrightarrow V(G)/I.\Lie \, E(G)$ de cette extension. 
\\

Considérons le diagramme 
$$
\xymatrix@R=6mm{
0 \ar[r] & \omega_G^*\otimes\widehat{\overline{K}} \ar[r] & 
V_p (G) \otimes_{\Qp} \widehat{\overline{K}} \ar[r]
\ar@{->>}[d]
 &
\omega_{G^D} \otimes\widehat{\overline{K}} \ar[r] \ar@{->>}[d]
& 0 \\
&& V_p (G)\otimes_F \widehat{\overline{K}}  \ar[r] &
\widetilde{\omega}_{G^D} \otimes \widehat{\overline{K}} 
}
$$
\begin{prop} \label{relhtO}
La suite 
$$
0 \ldrt \omega_G^*\otimes\widehat{\overline{K}} \ldrt  V_p
(G)\otimes_F \widehat{\overline{K}} \ldrt \widetilde{\omega}_{G^D}
\otimes \widehat{\overline{K}} \ldrt 0
$$
extraite du diagramme précédent est exacte.
\end{prop}
\dem
La surjectivité de l'application de droite est claire. De plus la
composée des deux applications est nulle. D'après la remarque B.9 de
l'appendice B de \cite{Cellulaire} 
$$
\dim \omega_G^*\otimes\widehat{\overline{K}} + \dim \widetilde{\omega}_{G^D}
\otimes \widehat{\overline{K}} = \dim V_p (G)\otimes_F \widehat{\overline{K}}
$$
Il suffit donc de montrer que l'application de gauche est injective
c'est à dire
$$
\omega_G^*\otimes\widehat{\overline{K}} \cap I.V_p(G)\otimes_{\Qp}
\widehat{\overline{K}} = 0
$$
Mais $I[\frac{1}{p}]$ est un produit de corps, il existe donc $e\in
I[\frac{1}{p}]$ tel que $e.I[\frac{1}{p}]= I[\frac{1}{p}]$. Mais
$I[\frac{1}{p}]. \omega_G^*\otimes\widehat{\overline{K}}=0$, d'où le
résultat.
\qed

\begin{rema} Comme dans la section précédente 
l'application $\Hom_\O (F/\O_F ,G) \ldrt \widetilde{\omega}_{G^D}$ se définit de la façon suivante en
termes de la $\O$-extension vectorielle universelle $E_\O (G)$ 
 : si $x:F/\O_F\ldrt
G_{O_K}$ alors il induit un morphisme 
$x_*:\O_{K}\ldrt \text{Lie} (E_\O (G))$
puisque $\O_{S/\Sigma}$ le faisceau structural du site cristallin
s'identifie au cristal algèbre de Lie de la $\O$-extension vectorielle universelle de $F/\O_F$.
\end{rema}

\section{Propriétés particulières de l'application de Hodge-Tate pour
  les groupes $p$-divisibles formels de dimension un}

On reprend les notations de la section précédente.

\subsection{Les périodes de Hodge-Tate vivent dans l'espace de Drinfeld}

Nous utiliserons la proposition suivante afin de définir 
l'isomorphisme au niveau des points.

\begin{prop}
Soit $G$ un groupe $p$-divisible 
sur $\O_K$. Considérons la suite de Hodge-Tate de $G^D$
$$
0 \ldrt \omega_{G^D}^*\otimes \mathcal{O}_{\widehat{\overline{K}}} (1)
\ldrt T_p (G)^* (1)\otimes_{\Zp} \mathcal{O}_{\widehat{\overline{K}}}  \ldrt 
\omega_G \otimes \mathcal{O}_{\widehat{\overline{K}}} \ldrt 0
$$ 
et plus particulièrement l'application 
$
\a_{G^D} : T_p (G)^* (1) \ldrt 
\omega_G \otimes \mathcal{O}_{\widehat{\overline{K}}}
$. Si $G$ est un groupe $p$-divisible formel c'est à dire  si sa fibre
spéciale $G_k$ ne
possède pas de partie étale alors l'application  $\a_{G^D}$ est injective.
\end{prop}
\dem 
Si $x:\Qp/\Zp\ldrt G_{/\O_{\overline{K}}}^D$ et $x^D :
G_{/\O_{\overline{K}}}\ldrt \mu_{p^\infty}$ alors, $x=0 \ssi
x^D=0$. De plus, $G$ étant formel $x^D=0$ ssi l'application tangente
associée est nulle, d'où le résultat.
\qed

\begin{coro}\label{DansOmega}
Soit  $G$ un $\O$-module $\pi$-divisible formel de dimension $1$ sur $\O_K$.  
Notons $\Omega \subset \mathbb{P} ( V_p (G)^*(1))$ l'espace de
Drinfeld au sens de Berkovich
obtenu en enlevant les hyperplans $F$-rationnels. Alors, avec les notations
précédentes 
$$
(V_p (G)^*(1)\otimes_{F}\widehat{\overline{K}} \twoheadrightarrow \omega_G \otimes 
 \widehat{\overline{K}}) \in \Omega ( \widehat{\overline{K}}) 
$$
\end{coro}

\subsection{Raffinement, d'après Faltings}

\begin{rema}
Soit $G$ un groupe p-divisible sur $\spf (\O_{K})$ et 
$$
0 \ldrt \omega_{G^D} \ldrt E (G) \ldrt G \ldrt 0
$$
son extension vectorielle universelle sur $\spf ((\O_{K})$. Alors, si 
\begin{eqnarray*}
E(G) (\O_{\widehat{\overline{K}}}) &=& \underset{n}{\limp} E(G) (\O_{\overline{K}}/ p^n \O_{\overline{K}}) \\
G(\O_{\widehat{\overline{K}}}) &=& \underset{n}{\limp} G(\O_{\overline{K}}/ p^n \O_{\overline{K}}) 
\end{eqnarray*}
il y a une suite exacte 
$$
0 \ldrt \omega_{G^D}\otimes \O_{\widehat{\overline{K}}} \ldrt E(G) (\O_{\widehat{\overline{K}}}) \ldrt G (\O_{\widehat{\overline{K}}}) \ldrt 0 
$$
Cela résulte aisément de la lissité des morphismes $E(G) \ldrt G$ 
réduits mod $\O_{\overline{K},\geq \l}$ pour tout $\l\in \Q_{>0}$ et de la surjectivité sur $\overline{k}$. 
Nous n'auront besoin de cela que pour $\Qp/\Zp$ pour lequel les points à valeurs dans une $\Zp$-algèbre de $E(\Qp/\Zp)$ sont $(R\oplus\Qp)/\Zp$, et cela est donc évident. 
\end{rema}

\begin{lemm}\label{fonda}
Soit $G$ un groupe p-divisible sur $\O_{K}$. Notons 
$$
\a_G : T_p (G) \ldrt \omega_{G^D}\otimes \O_{\widehat{\overline{K}}} 
$$
l'application de Hodge-Tate.
Soit $x=(x_n)_{n\geq 1} \in T_p (G)$ où $x_n\in G[p^n](\O_{\overline{K}})$. Soit $k$ un entier positif.  
 Alors, $\a_G (x)\in p^k \omega_{G^D}\otimes \O_{\widehat{\overline{K}}}$ ssi 
$x_k$ se relève en un point de $p^k$-torsion dans l'extension vectorielle universelle de $G$ i.e. un point de
$E(G)[p^n](\O_{\widehat{\overline{K}}})$.
\end{lemm} 
\dem
Considérons le diagramme 
$$
\xymatrix{
0 \ar[r] & \O_{\widehat{\overline{K}}} \ar[r] \ar[d]^{\times \a_G (x)} 
& (\O_{\widehat{\overline{K}}} \oplus \Qp)/\Zp \ar[r] \ar[d]^{E(x)}
 &  \Qp/\Zp \ar[r]\ar[d]^x & 0 \\
0 \ar[r] & \omega_{G^D} \otimes \O_{\widehat{\overline{K}}} \ar[r] & E(G) (\O_{\widehat{\overline{K}}}) \ar[r] & G( \O_{\widehat{\overline{K}}}) \ar[r] & 0
}
$$
L'élément $x_k\in G[p^k]$ se relève en un élément $\widetilde{x_k}= E(x) ([(0,p^{-k})]) \in E(G) (\O_{\widehat{\overline{K}}})$ et de plus 
$$p^k \widetilde{x_k}=\a_G (x)$$
Les relèvements de $x_k$ à $E(G)(\O_{\widehat{\overline{K}}})$ formant un $\omega_{G^D}\otimes 
\O_{\widehat{\overline{K}}}$-torseur le résultat s'en déduit.
\qed
\\

La proposition qui suit s'applique en particulier aux groupes $p$-divisibles formels
de dimension $1$ et est une légère amélioration d'un résultat de
Faltings. 

\begin{prop}[Faltings]
Soit $G$ un groupe $p$-divisible  formel sur $\O_K$. Considérons la suite de Hodge-Tate de $G^D$
$$
0 \ldrt \omega_{G^D}^*\otimes \mathcal{O}_{\widehat{\overline{K}}} (1)
\ldrt T_p (G)^* (1)\otimes  \mathcal{O}_{\widehat{\overline{K}}} \xrig{\; \a_{G^D}\otimes Id\;}
\omega_G \otimes \mathcal{O}_{\widehat{\overline{K}}} \ldrt 0
$$ 
 Supposons qu'il existe un groupe p-divisible isocline $H$ 
sur $k$ de hauteur $h$ tel que $\Z_{p^h}\hookrightarrow  \End
(H_{\bar{k}})$ 
et une isogénie $\rho : H\otimes_k \O_{K}/ p \O_{K} \ldrt G\otimes \O_{K}/p \O_{K}$ telle que $p^n \rho^{-1}$ soit une isogénie. Alors, si on note pour $M$ un ``réseau'' et $m\in M$ $\delta (m)= \sup \{k\;|\; p^{-k} m \in M\}$ on a 
$$
\forall x\in T_p(G)^*(1) \;\;\; \delta (x) \leq \delta (\a_{G^D} (x)) \leq \delta (x) + n +1 
$$
\end{prop}
\dem
Tout repose sur le lemme précédent. On peut supposer $k$ algébriquement clos. 
Soit $H_0$ un relèvement C.M. par $\Z_{p^h}$ de $H$ à $W(k)$.

Soit $x\in T_p (G^D) \setminus p T_p (G^D)$. Supposons que 
$\a_{G^D} (x)\in p^{n+1} \omega_G\otimes \O_{\widehat{\overline{K}}}$. Le lemme précédent implique que $x_{n+1}\in G^D [p^{n+1}](\O_{\overline{K}})$ se relève en un élément 
$$
\widetilde{x_{n+1}} \in E(G^D)[p^{n+1}]( \O_{\widehat{\overline{K}}})
$$
Il résulte de la nature cristalline de l'extension vectorielle
universelle que $\rho^D$ et $(p^n \rho^{-1})^D= p^n (\rho^D)^{-1}$ se relèvent en des morphismes 
$$
\xymatrix@C=2cm{
E(H_0^D)_{/\O_{K}} \ar@<1ex>[r]^{E(p^n (\rho^D)^{-1})} & \ar@<1ex>[l]^{E ( \rho^D)}  E (G^D)
}
$$
tels que $E(\rho^D) E( p^n (\rho^D )^{-1})=p^n$ et $ E( p^n (\rho^D
)^{-1}) E(\rho^D ) = p^n$. Considérons l'élément 
$y=E(\rho^D ) (\widetilde{x_{n+1}})\in E ( H_0^D) [p^{n+1}] ( \O_{\widehat{\overline{K}}})$.
Alors, $$y\neq 0$$ car sinon on aurait 
$$
p^n \widetilde{x_{n+1}}= E(p^n (\rho^D)^{-1}) (y) = 0 
$$
mais via $E(G^D)\drt G^D$, $p^{n} \widetilde{x_{n+1}}\mapsto x_1\in G[p]^D$ ce qui est en contradiction avec $x_1\neq 1$ puisque $x\notin pT_p (G^D)$. 

Donc, $\exists z\in E(H_0^D)[p] ( \O_{\widehat{\overline{K}}}),z \neq
0$. Soit $w\in H_0[p]^D (\O_{\overline{K}})$ son image. 
On a $w\neq 0$ car 
$$\ker ( E(H_0^D)[p] ( \O_{\widehat{\overline{K}}})
\ldrt H_0 [p]^D (  \O_{\widehat{\overline{K}}})) = (\omega_{H_0}\otimes
\O_{\widehat{\overline{K}}}) [p] =0
$$
 Soit $\iota : \Z_{p^h} \drt \End
(H_0^D)$ l'action C.M. obtenue par dualisation de Cartier de l'action
C.M. sur $H_0$. Le module de Tate $T_p (H_0^D)$ est $\Z_{p^h}$-libre de rang $1$. Donc, $\iota (\Z_{p^h}).w=H_0[p]^D(\O_{\overline{K}})$ qui est dans l'image de $E(H_0^D)[p]\drt H_0^D[p]$ puisque le morphisme $E(H_0^D)\drt H_0^D$ est $\Z_{p^h}$-équivariant. Par application du lemme \ref{fonda}
on en déduit que 
$$
\a_{H_0^D} (T_p (H_0^D)) \subset p \omega_{H_0} \otimes  \O_{\widehat{\overline{K}}}
$$
ce qui est en contradiction avec le fait que le conoyau du morphisme
de Hodge-Tate pour $H_0^D$ est annulé par un élément de valuation $\frac{1}{p-1}$.
\qed
\\

{\it Interprétation : } 
Dans les cas des groupes $p$-divisibles formels de dimension $1$ 
la proposition précédent dit plus précisément que si le point
$(G,\rho)$ est ``à distance $\leq n$'' du centre de l'espace de
Lubin-Tate alors ses périodes de Hodge-Tate dans l'espace de Drinfeld sont
``dans une boule de rayon $n+1$'' dans l'immeuble associé à l'espace
de Drinfeld.

\subsection{Formule exacte}

 Dans l'article \cite{Rami} nous donnons une formule exacte plus
 précise que la proposition précédente dans le cas des groupes formels de dimension $1$ 
pour le point associé dans
 l'espace de Drinfeld (cf. corollaire \ref{DansOmega}). Soit $|\mathcal{I}|$ la réalisation géométrique de  l'immeuble de Bruhat-Tits du groupe
 $\text{PGL}_n$. 
Il y a une application $s: \Omega (\widehat{\overline{K}})\ldrt |\mathcal{I}|$. 
 Cette formule exprime l'image par $s$ du point de Hodge-Tate 
en fonction  de la filtration de ramification sur
 le module de Tate du groupe $p$-divisible.

\section{Notations concernant les espaces de Lubin-Tate et de
  Drinfeld}
\label{kuipom}

Nous reprenons les notations du premier chapitre de \cite{Cellulaire}. 
Rappelons
que l'on note $\O=\O_F$. On note $L=W_\O(\Fqb)_\Q$ et $\s$ son Frobenius . 
 
Soit $\mathbb{H}$ un $\O$-module $\pi$-divisible formel de dimension $1$ et hauteur $n$ sur $\overline{\mathbb{F}}_p$
qui est défini sur $\mathbb{F}_q$. On a donc $\mathbb{H}=
\Hb^{(q)}$. On note alors
$\Pi=\text{Frob}_q\in \End (\Hb)$. Alors, 
$$
\O_D= \O_{F_n} [\Pi] \iso \End (\Hb )
$$
où l'on pose 
$$
\O_{F_n}= W_\O (\mathbb{F}_{q^n})
$$
Posons $\Gb=\Hb^n$ qui est un $\O$-module $\pi$-divisible formel de dimension $n$ et $\O$-hauteur $n^2$. Munissons le
d'une structure de $\O_D$-module formel spécial au sens de Drinfeld en posant 
$$
\iota: \O_D \ldrt \End ( \Gb) = M_n (\O_D)
$$
tel que 
$$
\forall x\in \O_{F_n} \;\; \iota ( x) = \text{diag} (x,x^\s,\dots, x^{\s^{n-1}} )
$$
et 
$$
\iota (\Pi) = \text{diag} ( \Pi,\dots,\Pi )
$$
et donc
$$
\forall x\in \O_D\;\; \iota (x)= \text{diag} ( x, \Pi x\Pi^{-1},\dots, \Pi^{n-1} x \Pi^{-(n-1)} )
$$
Tous les indices de $\Gb$ sont critiques et comme élément de $\text{PGL}_n (F) \bc \widehat{\Omega} ( \overline{\mathbb{F}_q})$ cela définit un point dans l'intersection de $n$ composantes irréductibles de la 
fibre spéciale. Il y a de plus une isogénie de $\O$-modules $\pi$-divisibles munis d'une action de $\O_D$
\begin{eqnarray}\label{ioD}
\text{Id} \oplus \Pi \oplus \dots \oplus \Pi^{n-1}: \Hb^n \ldrt \Gb
\end{eqnarray}
de degré $q^{\frac{1}{2} n (n-1)}$ où $\O_D$ agit diagonalement sur $\Hb^n$ (qui n'est pas spécial).
Il y a donc un isomorphisme d'isocristaux munis d'une action de $D$
$$
\mathbb{D} (\Gb)_\Q \iso \mathbb{D} (\Hb)_\Q^n
$$
où les modules de Dieudonné sont les modules 
de Dieudonné covariants relatifs à $\O$ (il s'agit de l'évaluation du
cristal algèbre de Lie de la $\O$-extension vectorielle universelle
sur $\O_L\twoheadrightarrow k$, cf. l'appendice B de \cite{Cellulaire}).
\\

{\it Convention :  Désormais on notera $\mathbb{D}$ le module de Dieudonné covariant
relatif à $\O$ d'un $\O$-module $\pi$-divisible noté 
$\mathbb{D}_\O$ dans l'appendice B de \cite{Cellulaire} et par cristal de 
Messing on entendra le cristal de Messing $\O$-extension vectorielle
universelle  défini dans l'appendice B de 
\cite{Cellulaire}.}

\subsection{Modules de Dieudonné}

{\it Nous n'utiliserons pas cette sous-section plus tard.}
\\
Il y a une identification 
$$
\mathbb{D} (\Hb) = \O_D\otimes_{\O_{F_n}} \O_L 
$$
muni de $\iota : \O_D \ldrt \End ( \Hb)$ tels que si $V$ désigne le Verschiebung
\begin{eqnarray*}
\forall d\otimes \l \in  \mathbb{D} (\Hb) \;\;\, V ( d\otimes \l) & = & d \Pi \otimes \l^{\s^{-1}} \\
\forall d'\in \O_D\;\;\; \iota (d') ( d\otimes \l) & =& d' d \otimes \l 
\end{eqnarray*}
Alors, 
$$
\mathbb{D} (\Gb) = \O_D\otimes_{\O_F} \O_L 
$$
muni de son action de $\iota$ de $\O_D$ et 
\begin{eqnarray*}
\forall d\otimes \l \in \mathbb{D} (\Gb) \;\;\; V ( d\otimes \l) & = & d \Pi \otimes \l^{\s^{-1}} \\
\forall d'\in \O_D\;\;\; \iota (d') ( d\otimes \l) & =& d' d \otimes \l 
\end{eqnarray*}
Dès lors on a un isomorphisme de $\O_L$-modules
\begin{eqnarray*}
\O_D\otimes_{\O_F} \O_L & \xrig{\;\sim \;} & \prod_{k\in \Z/n\Z}
\O_D\otimes_{\O_{F_n},\s^{k}} \O_L \\
d\otimes \l & \longmapsto & ( d \otimes \l)_{k\in \Z/n\Z}
\end{eqnarray*}
 Et dans ces coordonnées l'identification $\mathbb{D} (\mathbb{H})^n = \mathbb{D} (\mathbb{G})$ est donnée par 
\begin{eqnarray*}
(\O_D\otimes_{\O_{F_n}} \O_L )^n &\ldrt & \prod_{k\in \Z/n\Z} \O_D
  \otimes_{\O_{F_n},\s^{k}} \O_L \\
(x_k\otimes\l_k)_{0\leq k \leq n-1} &\longmapsto & (\Pi^{-k} x_k\Pi^k\otimes \l_k)_{k\in\Z/n\Z} 
\end{eqnarray*}
et l'isogénie (\ref{ioD})  par
\begin{eqnarray*}
(\O_D\otimes_{\O_{F_n}} \O_L )^n &\ldrt & \prod_{k\in \Z/n\Z} \O_D
  \otimes_{\O_{F_n},\s^{k}} \O_L \\
(x_k\otimes\l_k)_k &\longmapsto & (x_k\Pi^k\otimes \l_k)_k 
\end{eqnarray*}

\subsection{Notations concernant les espaces de Lubin-Tate}

Rappelons que l'on note $\breve{F} =\widehat{F^{nr}}$
le complété de l'extension maximale non-ramifiée de $F$ dans une
clôture algébrique de celui-ci. 
 On fixe un
isomorphisme entre $\overline{\mathbb{F}}_q$ et le corps résiduel de
$\breve{F}$. Celui-ci
induit un isomorphisme $L\simeq \breve{F}$. Néanmoins on ne confondra 
pas toujours ces deux corps, l'un, $L$, étant un corps ``abstrait'',
l'autre, $\breve{F}$, un corps plongé (un corps reflex).

\begin{defi}
Soit $K$ un corps valué complet (pour une valuation de rang $1$ c'est
à dire à valeurs dans $\R$) extension de $\breve{F}$. 
On appelle point à valeurs dans $K$ de la tour de Lubin-Tate 
un triplet $(H,\rho,\eta)$ à isomorphisme près où 
\begin{itemize}
\item $H$ est
un $\O$-module $\pi$-divisible formel de dimension $1$  sur $\mathcal{O}_{K}$
\item $\rho$ est une quasi-isogénie $\mathbb{H}\otimes_{\overline{\mathbb{F}}_q} \O_{K}/ p \O_{K} \drt H\otimes   \O_{K}/ p \O_{K}$ 
\item $\eta : \O^n \iso T_p (H)$ est un isomorphisme de
   $\Gal(\overline{K}|K)$-modules, le membre de gauche étant muni de
   l'action triviale de Galois. 
\end{itemize}
On note $\M^{\LT}_\infty (K)$ cet ensemble qui est muni d'une action
de $D^\times \times \GL_n (\O_F)$, où $D^\times $ agit à gauche par
$d.\rho=\rho \circ d^{-1}$ et $\GL_n (\O_F)$ à droite par $\eta.g=\eta\circ g$, 
 et d'une donnée de descente $\a$ vers $F$.
\end{defi}

On renvoie à la section 1.4 de \cite{Cellulaire} pour la définition de
la donnée de descente $\a$.
On verra bientôt que cet ensemble est en fait muni d'une action de
$\GL_n (F)$. 

\begin{rema}
 Dans la définition précédente, la définition de $\eta$ peut être
 remplacée par : $\eta$ est un isomorphisme 
$$
\eta:\O^n\iso \Hom_\O ( F/\O, H)
$$
\end{rema}

\begin{defi}\label{ruygop}
On note $\M^{\LT} (K)/\sim$ l'ensemble des couples $(H,\rho)$ comme
précédemment à isogénie déformant un élément de $F^\times$  près i.e.
$(H,\rho)\sim (H',\rho') \ssi \exists f : H\ldrt H'$ une isogénie
et un $x\in F^\times$ tels que le diagramme suivant commute 
$$
\xymatrix{
H \text{ mod } \pi \ar[r]^{f \text{ mod } \pi} & H'\text{ mod } \pi \\
\mathbb{H}\otimes \O_K/\pi\O_K \ar[r]^{x}\ar[u]^\rho & \mathbb{H}\otimes
\O_K/\pi\O_K \ar[u]^{\rho'}
}
$$
Cet ensemble est muni d'une action de $D^\times$.

Une définition équivalente consiste à remplacer à isogénie déformant
un élément de $F^\times$ près par  à isogénie déformant
une puissance de $\pi$ près.

Une troisième définition équivalente consiste à dire que $(H,\rho)$
est équivalent à $(H',\rho')$ ssi il existe une quasi-isogénie
$f:H\ldrt H'$ sur $\spec (\O_K)$ telle que $f\rho=\rho'$ i.e. ssi la
quasi-isogénie $\rho'\rho^{-1}$ sur $\spec (\O_K/p\O_K)$ se relève en
une quasi-isogénie sur $\spec(\O_K)$.
\end{defi}

\subsection{Notations concernant les espaces de Drinfeld}

\begin{defi}
Soit $K$ un corps valué complet extension de $\breve{F}$. On appelle point à valeur dans $K$ de la tour de Drinfeld un triplet $(G,\rho,\eta)$ à isomorphisme près
où
\begin{itemize}
\item $G$ est un $\O_D$-module $\pi$-divisible formel spécial sur $\spf (\O_{K})$
\item $\rho$ est une quasi-isogénie $\O_D$-équivariante $\mathbb{G} \otimes_{\overline{\mathbb{F}}_q} \O_{K}/p\O_{K}\ldrt 
G\otimes_{\O_{K}}\O_{K}/p\O_{K} $
\item $\eta: \O_D\iso T_p (G)$ est un isomorphisme de $\O_D$-modules
  galoisiens.
\end{itemize}
On note $\M^\D_\infty (K)$ cet ensemble 
 muni de son action de $\GL_n (F)\times D^\times$ où 
$\GL_n (F)= \End_{D-\text{éq.}}(\mathbb{G})_\Q^\times$ agit à gauche 
via $g.\rho=\rho\circ g^{-1}$
et $\O_D^\times$ agit à droite via $\eta.d = \eta \circ d^{-1}$ où $
d^{-1}: \O_D\iso \O_D$ est $d'\mapsto d' d^{-1}$.  L'action de $\Pi\in
D^\times$ est définie par $\Pi. (G,\rho,\eta) = (G/G[\Pi], \ph \circ
\rho, \eta.\Pi)$ où $\ph : G\twoheadrightarrow G/G[\Pi]$ et 
$$
\xymatrix@C=15mm{
\O_D \ar[d]^{.\Pi} \ar[r]^\eta_\sim & T_p (G) \ar@{^(->}[d]^{\ph_*} \\
\O_D \ar[r]^{\eta.\Pi}_\sim & T_p ( G/G[\Pi]) 
}
$$
\end{defi}

\begin{rema}
On renvoi à la section 1 de \cite{iso4} pour la définition générale
de l'action d'un élément de $D^\times$. On notera que $\{(x,x)\;|\;
x\in F^\times \}\subset \GL_n (F)\times D^\times$ agit trivialement. 
\end{rema}

\begin{rema}
 Dans la définition précédente, la définition de $\eta$ peut être
 remplacée par :  on se donne un élément 
$$
\eta (1)\in \Hom_\O (F/\O, G)\setminus \Pi. \Hom_\O (F/\O, G)
$$
\end{rema}

\begin{rema} 
Le groupe associé à l'espace de Rapoport-Zink précédent et noté $G$ dans \cite{RZ}  est dans notre cas $(D^{opp})^\times$. Dans la définition précédente on l'a identifié à $D^\times$ via $d\mapsto d^{-1}$ ce qui explique que l'on pose $d.\eta = \eta \circ d^{-1}$ alors que la définition usuelle des correspondances de Hecke est pour $g\in G(\Qp)$, $\eta \mapsto \eta \circ g$.
\end{rema}

\begin{rema}
L'identification $\GL_n (F)= \End_{D-\text{éq}} (\mathbb{G})_\Q^\times$ se déduit, avec les notations des sections qui suivent, de l'identification
$$
 \End_{D-\text{éq}} (\mathbb{G})_\Q =  \End_{D-\text{éq}}\left
   (\mathbb{D}(\mathbb{G})_\Q, V \right ) = \End_\O \left ( \mathbb{D}(\mathbb{G})_{\Q,0}^{V^{-1}\Pi}\right )
$$
et d'un choix d'une base de $\mathbb{D}(\mathbb{G})_{\Q,0}^{V^{-1}\Pi}$ (cf. les section suivantes).
\end{rema}

\begin{defi}
Comme dans la définition \ref{ruygop} 
on note $\M^{\D} (K)/\sim$ l'ensemble des couples $(G,\rho)$  à isogénie déformant un élément de $F^\times$ près. Cet ensemble
est muni d'une action de $\GL_n (F)$.
\end{defi}

\subsection{Quelques rappels sur Drinfeld classique}

Notons  
$$
N^{\D}= \mathbb{D} (\Gb)_\Q = D\otimes_F L
$$
un isocristal relativement à l'extension $L|F$. 
On a 
\begin{eqnarray*}
N^\D & =& \bigoplus_{i\in \Z/n \Z} N^\D_i \;\text{ où }  N^\D_i =\{n\in N\; |\; \forall x\in F_n \;
\iota (x).n= x^{\s^{-i}} n \; \} 
\end{eqnarray*}
où relativement à cette graduation $\deg V= \deg \Pi =+1$. 
On a de plus 
\begin{eqnarray*}
D\otimes_F L  &=& \bigoplus_{a\in \Z/n\Z} D\otimes_{F_n,\s^{-a}} L \\
&=& \bigoplus_{a,b\in \Z/n\Z}
L. \underbrace{\iota(\Pi^b)(e_a)}_{e_{a,b}} 
\end{eqnarray*}
où $e_a =1\otimes 1\in D\otimes_{F_n,\s^{-1}}L$ et $e_{a,b} =\Pi^b\otimes 1$. 
De plus, avec ces coordonnées
$$
V(e_{a,b})= \iota(\Pi)(e_{a,b}) = \pi^{\delta_{b,n-1}} e_{a,b+1}
$$
et 
$$
 N^\D_i = \bigoplus_{a+b =i} L.e_{a,b}
$$
 L'opérateur
$$
V^{-1} \Pi : N_0^\D\ldrt N^\D_0
$$
est de pente $0$ et fait de $(N_0^\D,V^{-1} \Pi)$ un isocristal unité :
\begin{eqnarray*}
(N^\D_0, V^{-1} \Pi)\simeq ((N^\D_0)^{V^{-1} \Pi}\otimes_F L, Id\otimes \s) \\
(N^\D)^{V^{-1} \Pi}_0 = \bigoplus_{a+b = 0} F.e_{a,b} \simeq F^n
\end{eqnarray*}

Rappelons maintenant qu'il y a des bijections 
\begin{eqnarray*}
\left \{
{\text{Sous isocristaux D-stables de} \atop (N^\D,\ph)} \right \} &\iso &
\left \{ \text{Sous } F\text{-ev. de} \atop
(N^\D)^{V^{-1} \Pi} \right \} \\
N & \longmapsto & N_0^{V^{-1} \Pi}
\end{eqnarray*}
où $\ph$ désigne le Frobenius, 
et pour $K|L$ une extension valuée comme dans le premier chapitre 
\begin{eqnarray*}
\left \{ \text{ Filtrations D-stables de codimension } n \atop \text{dans } N^\D\otimes_L K
\right \} &\iso  & \left \{ \text{Filtrations de codimension } 1 \text{ dans}\atop 
(N_0^\D)^{V^{-1} \Pi} \otimes_F K \right \} \\
\Fil=\oplus_{i\in \Z/n \Z} \Fil_i &\longmapsto & \Fil_0
\end{eqnarray*}
Alors, $(N^\D,\ph,\Fil)$ est faiblement admissible ssi pour tout sous-isocristal 
$N\subset N^\D$ on a 
$$
t_H ( N, N\otimes_L K\cap \Fil) \leq t_N (N,\ph)
$$
où $t_H$, resp. $t_N$, désigne le point terminal du polygone de Hodge,
resp. Newton.
L'existence de la filtration de Harder-Narashiman et sa canonicité sur $(N^\D,\ph,\Fil)$ impliquent que celle-ci est $D$-stable si $\Fil$ l'est et qu'il suffit donc dans ce cas là de tester la condition 
d'admissibilité faible sur les sous-isocristaux $D$-stables
(\cite{RZ} chapitre 1). Restreignons-nous aux filtrations $\Fil$ $D$-stables de codimension $n$. 
 Via les deux bijections ci dessus 
on trouve aisément :
$$
\left \{ 
\text{Fil} \subset N^\D\otimes_L K \text{ telles que } (N^\D,\ph,\Fil ) \atop
  \text{ est faiblement admissible }  \right \}
\iso
\left \{ \Fil'\subset (N^\D_0)^{V^{-1} \Pi}\otimes K \;   \text{ telles que } \atop
\forall N'\subset (N^\D_0)^{V^{-1} \Pi} \text{ de dim.} 1 \; \Fil'\cap (N'\otimes K) = (0)
\right \}
$$
la condition portant donc sur les sous-isocristaux $D$-stables de dimension $n$ dans
$(N^\D,\ph)$ qui sont paramétrés par $\check{\mathbb{P}}\left ( (  N^\D_0)^{V^{-1} \Pi} \right )(F)$ où pour $E$ un $F$-e.v. on note $\mathbb{P} (E)$ l'espace des filtrations de codimension $1$ de $E$ et $\check{\mathbb{P}} (E) = \mathbb{P} \left ( \check{E} \right )$ l'espace des droites de $E$. La relation d'incidence entre $\mathbb{P} (E)$ et  $\check{\mathbb{P}} (E)$ est donnée par : pour $H\in \check{\mathbb{P}} (E)$, $H$ définit une droite $D$ dans $E$ et donc un hyperplan  encore noté $H$ dans  $\mathbb{P} (E)$ formé des filtrations contenant $D$. Ainsi $ \check{\mathbb{P}} (E)$ paramètre les hyperplans dans  $\mathbb{P} (E)$. 
 On a 
donc 
$$ \left 
\{ (N^\D,\ph,\Fil) \text{ faiblement admissible} \atop \Fil \; D\text{-stable}
\right \} \iso
{\mathbb{P}}\left ( N^\D_0\right )(K)\setminus \bigcup_{H\in \check{
\mathbb{P}}\left  ((N^\D_0)^{V^{-1} \Pi}\right )(F)}  H(K)  = \Omega (K)
$$
où $\Omega$ désigne l'espace de Drinfeld vu comme espace de
Berkovich. 

Dit d'une autre façon, $\Omega (K)$ est l'ensemble des filtrations  $\Fil$
de codimension $1$ dans $N^\D_0\otimes_L K$ telles que l'application suivante soit injective
$$
(N^\D_0)^{V^{-1} \Pi} \hookrightarrow  N^\D_0\otimes_L K/\Fil 
$$

\section{Description de $\M^{\LT} (K)/\sim$ en termes de modules filtrés}

Il s'agit ici d'étendre la théorie de Fontaine des modules faiblement admissibles 
pour les groupes de dimension $1$ 
à des corps $K$ non forcément de valuation
discrète.

Soit un couple  $(H,\rho)$ où $ H$ est un $\O$-module $\pi$-divisible  sur $\spf (\O_K)$ et $\rho$ une 
rigidification avec $\mathbb{H}$ modulo $p$. L'algèbre de Lie de la $\O$-extension vectorielle universelle
de $H$, $M$, est munie d'une filtration $\Fil \subset M$ telle que
$M/\Fil  \simeq \omega_H^*$. La nature cristalline de la
$\O$-extension vectorielle universelle induit un isomorphisme 
$$
\rho_* : \mathbb{D} (\Hb)_\Q \otimes_{L} K \iso M\unp
$$
et $\Fil_H= \rho_*^{-1} (\Fil\unp) \subset\mathbb{D} (\Hb)_\Q \otimes K$ est une filtration de codimension $1$.

\begin{prop}\label{four_D}
L'application précédente $(H,\rho_H)\mapsto \Fil_H$ induit une bijection  $D^\times$-équivariante
$$
 \M^{\LT} (K)/\sim \iso \P \left ( \mathbb{D} (\Hb)_\Q \right ) (K)
 $$
\end{prop}
\dem
C'est une conséquence de l'existence du domaine fondamental de Gross Hopkins et du fait que la 
restriction du morphisme des périodes à ce domaine est un isomorphisme
d'espaces rigides (cf. \cite{HopkinsGross} et \cite{Cellulaire}).

Plus précisément, étant donné une filtration $\Fil \in \mathbb{P}(K)$ il existe un entier $i$ tel que 
$\Pi^i.\Fil$ soit dans l'image du domaine fondamental de Gross-Hopkins. Il existe donc un couple $(H,\rho)\in \M^{\LT} (K)$
induisant  $\Pi^i.\Fil$. Alors, $\Pi^{-i}.(H,\rho)$ induit $\Fil$, d'où la surjectivité de l'application.

Pour l'injectivité, si $(H,\rho)$ et $(H',\rho')$ induisent la même filtration alors la quasi-isogénie 
$$
\rho'\circ \rho^{-1} : H\otimes {\O_K/p\O_K} \ldrt H'\otimes \O_K/p\O_K
$$
est telle que l'application induite 
$$
(\rho'\circ \rho^{-1})_* : \text{Lie} (E(H))\ldrt \text{Lie} ( E ( H'))
$$
vérifie $$\exists a\in \N\; (\rho'\circ \rho^{-1})_* (\omega_{H^D}) \subset p^{-a} \omega_{{H'}^D}$$
et donc d'après la théorie de Messing l'isogénie 
$$
p^a \rho'\circ \rho^{-1} : (H,\rho) \ldrt (H',\rho') \text{  mod } p
$$
sur $\spec(\O_K/p\O_K)$ se relève en
 une isogénie : $(H,\rho)\sim (H',\rho')$.
\qed
 
\begin{rema}
  Il se peut qu'en général pour des espaces de périodes plus généraux l'application des périodes ne soit pas surjectivité sur les $K$-points  ! i.e. qu'un module faiblement admissible ne provienne pas forcément d'un groupe p-divisible  sur $K$ quelconque. 
Par contre le calcul des fibres de l'application des périodes (la
partie injectivité dans la démonstration précédente) reste valable en général.
\end{rema}

\section{Description de $\M^\D (K)/\sim$ en termes de module filtré}

Soit $(G,\rho)$ où $G$ est un $\O_D$-module formel spécial sur $\spf (\O_K)$ et $\rho$ une rigidification avec $\mathbb{G}$. Soit $M$ l'extension vectorielle universelle de $G$ qui est
filtrée via $\Fil\subset M$ où $M/\Fil \simeq \omega_G^*$. Il y a une décomposition sous l'action de $F_n$ 
$$
M=\bigoplus_{i\in \Z/n \N} M_i \;\text{ et } \; \Fil= \bigoplus_{i\in\Z/ n\Z} \Fil_i 
$$
d'où un élément $\Fil_{G,0} = \rho_*^{-1} (\Fil _0) \subset \mathbb{D} (\Gb)_{\Q,0} \otimes K, \;\Fil_{G,0}\in \Omega (K)$.

\begin{prop}\label{DrP}
L'application $(G,\rho)\mapsto \Fil_{G,0}$ induit une bijection $\GL_n (F)$-équivariante 
 entre $\M^\D (K)/\sim$ et $\Omega (K)$ avec
  $$\Omega (K) =
{\mathbb{P}} \left ( \mathbb{D} (\Gb)_{\Q,0} \right )(K)\setminus \bigcup_{H\in \check{
\mathbb{P}} ((\mathbb{D}(\Gb)_{\Q,0})^{V=\Pi})(F)}  H(K)  
$$
\end{prop}
\dem Cela résulte de ce que l'application des périodes est un isomorphisme rigide pour les
espaces de Drinfeld, de l'égalité $\Omega (K)=\widehat{\Omega} (\O_K)$ où $\widehat{\Omega}$ est le schéma formel de Drinfeld, du théorème 
de Drinfeld et du fait que les deux applications de périodes, celle
définie par Drinfeld et celle définie dans \cite{RZ},  coïncident.
\qed

\section{Prolongement des isogénies}

\subsection{Prolongement}

Soit $H$ un groupe $p$-divisible sur $\spec (\O_K )$. Soit $\mathcal{C}$ l'ensemble des classes d'isomorphismes de
couples $(H',f)$ où $H'$ est un groupe
$p$-divisibles sur $\O_K$ et $f:H'\drt H$ une quasi-isogénie sur $\spec (\O_K)$. 
 Pour un tel couple $T_p (f): T_p (H')\hookrightarrow V_p (H)$. 

\begin{lemm}\label{adsch}
L'application 
\begin{eqnarray*}
\mathcal{C} & \drt & \{ \text{ réseaux } \Gal (\overline{K}|K)-\text{stables dans } V_p (H)\; \} \\
(H',f) & \mapsto & Im \;T_p (f)
\end{eqnarray*}
est une bijection.
\end{lemm}
\dem
C'est une conséquence de ce que pour $\mathcal{G}$ un groupe fini
localement libre sur $\O_K$ et $I\subset \mathcal{G}_\eta$ un
sous-groupe fini localement libre il existe un unique prolongement
$\mathcal{I}\hookrightarrow \mathcal{G}$ de l'inclusion $I\subset
\mathcal{G}_\eta$ où $\mathcal{I}$ est un sous-groupe fini localement
libre. Le groupe $\mathcal{I}$ est obtenu par adhérence schématique
(cf. le chapitre 2 de \cite{Ray1} pour le cas de valuation discrète 
et le lemme qui suit en général).
\qed

\begin{lemm}
Soit $E\subset K^n$ un sous-$K$-espace vectoriel. Alors, $E\cap \O_K^n$ est un $\O_K$-module libre de rang fini facteur direct.
\end{lemm}

\subsection{Définition de l'action de $\GL_n (F)$ sur $\M^{\LT}_\infty (K)$}

On a définit une action de $\GL_n (\O_F)$ sur $\M^{\LT}_\infty (K)$. 
Utilisons le lemme précédent pour étendre cette action à $\GL_n (F)$. 
Soient $[(H,\rho,\eta)]\in \M_\infty^{\LT} (K)$ et $g\in \GL_n
(F)$ (les crochets signifient que l'on prend une classe d'isomorphisme
de triplets). D'après le lemme précédent au réseau Galois-stable $\eta (g.
\O_F^n)$ dans $V_p (H)$ correspond une quasi-isogénie $\ph :H\ldrt
H'$ sur $\spec (\O_K)$ telle que $\ph_*^{-1} ( T_p (H')) = \eta (
g.\O_F^n)$
où $\ph_* : V_p (H)\iso V_p (H')$. 
 On pose
alors $$g.[(H,\rho,\eta)] = [(H',(\ph \text{ mod } p)\circ \rho,\ph_*\circ
\eta \circ
g)$$

\section{Description de $\M^{\LT}_\infty (K)$ en termes de
  modules filtrés rigidifiés}

\begin{defi}
Soit $[(H,\rho)]\in \M^{\LT}(K)/\sim$. 
Une rigidification du module de Tate de la classe d'isogénie $[(H,\rho)]$ est un isomorphisme de modules galoisiens
$$
\eta : F^n\iso V_p (H)
$$
qui induit donc naturellement $\forall (H',\rho')\sim (H,\rho)$ un
isomorphisme
$$
F^n\iso V_p (H')
$$
via l'isomorphisme canonique $V_p (H)\iso V_p (H')$ induit par le
relèvement de la quasi-isogénie $\rho'\rho^{-1}$ sur $\spec (\O_K)$. 
\end{defi}

\begin{prop}\label{peftu}
L'application naturelle 
$$
\M_\infty^{\LT} (K)\ldrt \{ \left ([(H,\rho)],\eta \right )\; |\; [(H,\rho)]\in
\M^{\LT} (K)/\sim\text{ et } \eta \text{ une rigidification} \;\}
$$
est une bijection $\GL_n (F)\times D^\times$-équivariante.
\end{prop}
\dem
Décrivons l'inverse de cette application. Soit $([(H,\rho)],\eta)$
dans le membre de droite. D'après le lemme \ref{adsch} il existe un
unique $(H',\rho')\sim (H,\rho)$ tel que 
$$
\eta : \O_F^n\iso T_p (H')
$$
On associe alors à  $([(H,\rho)],\eta)$ le triplet
$(H',\rho',\eta')$. On vérifie aisément que ces deux applications sont
inverses l'une de l'autre.
\qed

Ainsi les fibres de l'application $\M_\infty^{\LT} (K) \ldrt \M^{\LT}
(K)/\sim$ sont des $\GL_n (F)$-torseurs. 
Soit $U\subset \GL_n (\O_F)$ un sous-groupe compact-ouvert et
$\M_U^{\LT} (K)$ l'ensemble des $K$-points de l'espace de Lubin-Tate en niveau $U$ (que nous
n'avons pas défini dans cet article).  
En fait on a la décomposition
plus précise d'extensions ``Galoisiennes'' 
$$
\xymatrix{
    \ar@{-}@/_4pc/[dd]_{GL_n (F)} \M_\infty^{\LT} (K) \ar[d]  \ar@{-}@/^2pc/[d]^{U} \\
  \M_U^{\LT} (K) \ar[d] \\
  \M^{\LT} (K)/\sim \ar[r]^{\sim} & \mathbb{P} ( \mathbb{D} (\Hb)) (K) 
}
$$
et on peut donc ainsi dire que l'espace des périodes de Gross-Hopkins
est le quotient de l'espace de Lubin-Tate en niveau infini par $\GL_n
(F)$.

\begin{theo}\label{derhj}
Il y a une bijection $\GL_n (F)\times D^\times$-équivariante entre
l'ensemble $\M_\infty^{\LT} (K)$ et l'ensemble des couples
$(\Fil,\zeta)$ où
$$
\Fil\in\mathbb{P} \left (\mathbb{D} (\mathbb{H} \right )_\Q) (K)
$$
$\zeta= (\zeta_1,\dots,\zeta_n)$ où $$\forall i\; \zeta_i\in\Hom_\O
(F/\O_F,\mathbb{H}\otimes_{\overline{\mathbb{F}}_p} \O_K/p\O_K)\unp $$
sont linéairement indépendants sur $F$ i.e. induisent une injection 
$F^n\hookrightarrow \Hom_\O (F/\O_F,\mathbb{H}\otimes \O_K/p\O_K)\unp$,
et $\forall i\;$ le morphisme induit sur l'évaluation du cristal de
Messing sur l'épaississement $\O_K\twoheadrightarrow \O_K/p\O_K$,
$\zeta_{i*} : K\drt \mathbb{D} (\mathbb{H})\otimes K$, vérifie
  $\zeta_{i*}(1)\in \Fil$.
\end{theo}
\dem
Il s'agit d'une conséquence de la proposition précédente couplée 
à la proposition \ref{four_D} et au
critère de relèvement de Messing vis à vis de l'idéal à puissances
divisées engendré par $p$. Dans cet énoncé, si $(e_i)_{1\leq i \leq
  n}$ est la base canonique de $F^n$ on a posé 
$$
\zeta_i = \rho^{-1} \circ \left ( \eta (e_i) \text{ mod } p \right)
$$
\qed

\subsection{Description de $\M_\infty^{\LT} (K)$ en termes d'algèbre
  linéaire}\label{bujiko}

Cette sous-section ne sera pas nécessaire dans la suite mais donne une
description purement en termes d'algèbre linéaire de $\M_\infty^{\LT}
(K)$. 

\begin{theo}\label{deall}
Supposons $F=\Qp$. 
Il y a une bijection $\GL_n (F)\times D^\times$-équivariante entre
l'ensemble $\M_\infty^{\LT} (K)$ et l'ensemble des couples
$(\Fil,\xi)$ tels que 
$$
\Fil\in\mathbb{P} (\mathbb{D} (\mathbb{H})_\Q) (K)
$$ et 
$$
\xi : F^n\hookrightarrow \Fil\left ( \mathbb{D}(\mathbb{H})_\Q\otimes_L
B^+_{cris} (\O_K)\right )^{\ph = p}
$$
(qui est alors automatiquement un isomorphisme)
\end{theo}
\dem 
Il s'agit d'une conséquence de la proposition \ref{peftu} couplée à la
proposition 
\ref{four_D} et au théorème \ref{crop}. 
\qed 

\section{Description de $\M_\infty^{\D} (K)$ en termes de modules
  filtrés rigidifiés}

\begin{defi}
Soit $[(G,\rho)]\in\M^{\D} (K)/\sim$. Une rigidification de la classe
d'isogénie $[(G,\rho)]$ est un isomorphisme de $D$-modules galoisiens 
$$
\eta : D\iso V_p(G)
$$
qui induit donc naturellement $\forall (G',\rho')\sim (G,\rho)$ un
isomorphisme 
$$
D\iso V_p (G')
$$
via l'identification $V_p (G)\iso V_p (G')$ induite par le relèvement
de $\rho'\rho^{-1}$.
\end{defi}

\begin{rema}\label{huyto}
Dans la proposition précédente la donnée de $\eta$ est équivalente à
celle de $\zeta = \eta (1) \in \Hom_\O (F/\O_F, G) \unp\setminus \{ 0 \}$.
\end{rema}

\begin{prop}
L'application naturelle
$$
\M_\infty^{\D} (K) \ldrt \{\; \left ([(G,\rho)],\eta \right )\; |\;
  [(G,\rho)]\in\M^{\D}(K)/\sim \text{ et } \eta \text{ une rigidification}\;\}
$$
est une bijection $\GL_n (F)\times D^\times$-équivariante.
\end{prop}
\dem
Elle est identique à celle de la proposition \ref{peftu} en utilisant
le lemme \ref{adsch}.
\qed
\\

Comme pour l'espace de Lubin-Tate il y a une description pour
$U\subset \O_D^\times$ un sous-groupe ouvert
$$
\xymatrix{
    \ar@{-}@/_4pc/[dd]_{D^\times} \M_\infty^{\D} (K) \ar[d]  \ar@{-}@/^2pc/[d]^{U} \\
  \M_U^{\D} (K) \ar[d] \\
  \M^{\D} (K)/\sim \ar[r]^(.6){\sim} & \Omega (K) 
}
$$

\begin{theo}\label{descDri}
Il y a une bijection $\GL_n (F)\times D^\times$-équivariante entre
l'ensemble $\M^{\D}_\infty (K)$ et l'ensemble des couples
$(\Fil,\zeta)$ où
$$
\Fil\in \Omega (K)\subset \mathbb{P}\left
  (\mathbb{D}(\mathbb{G})_{\Q,0}\right )(K)
$$
et 
$$
\zeta \in \Hom_\O (F/\O_F, \mathbb{G}\otimes \O_K/p\O_K)\unp \setminus
\{ 0\}
$$
est tel que si $\zeta_* : K\drt \mathbb{D} (\mathbb{G})_{\Q}\otimes K$
est l'application induite sur l'évaluation du cristal de Messing alors 
si 
$$
\zeta_* (1) = \oplus_i \zeta_*(1)_i\in \bigoplus_{i\in \Z/n\Z}
\mathbb{D}(\mathbb{G})_{\Q,i}\otimes K
$$ 
on a 
$$
\forall i\;\; \Pi^{-i} \zeta_* (1)_i \in \Fil
$$
\end{theo}
\dem
C'est une conséquence de la proposition précédente couplée à la 
remarque \ref{huyto} et au
théorème de relèvement de Messing. Dans l'énoncé on a posé 
$$
\zeta = \rho^{-1} \circ \left ( \eta (1) \text{ mod }  p \right )
$$
\qed

\subsection{Description de $\M_\infty^{\D} (K)$ en termes d'algèbre linéaire}

Comme dans la section \ref{bujiko} cette section ne servira pas dans
la suite.

\begin{theo}
Supposons $F=\Qp$. 
Il y a une bijection $\GL_n (F)\times D^\times$-équivariante entre
l'ensemble $\M^{\D}_\infty (K)$ et l'ensemble des couples
$(\Fil,\xi)$ où
$$
\Fil\in \Omega (K)\subset \mathbb{P}(\mathbb{D}(\mathbb{G})_{\Q,0})(K)
$$
et 
$$
\xi \in  \Fil(\mathbb{D}(\mathbb{G})_{\Q,0}\otimes
B^+_{cris}(\O_K))^{(V^{-1}\Pi\otimes \ph)^n=p}\setminus \{ 0 \} 
$$
\end{theo}

\section{La bijection au niveau des points}\label{coeulop}

On fixe
$$
\Delta : \Hb^n \ldrt\Gb
$$
une quasi-isogénie compatible à l'action de $\O_D$ comme par exemple
celle définie dans la section \ref{kuipom}.

On fixe  un isomorphisme 
\begin{eqnarray}\label{out1}
\mathbb{D}
(\mathbb{H})_{\Q,0}^{V=\Pi} \simeq F
\end{eqnarray}
 qui induit via $\Delta$ un
isomorphisme 
\begin{eqnarray*}\label{out2}
\mathbb{D} (\mathbb{G})_{\Q,0}^{V=\Pi} \simeq F^n
\end{eqnarray*}
et induit donc 
\begin{eqnarray*}\label{out3}
\mathbb{D} (\mathbb{G})_{\Q,0} \simeq L^n
\end{eqnarray*}
De plus,  l'isomorphisme (\ref{out1}) induit également un isomorphisme 
$\mathbb{D}(\mathbb{H})_{\Q,0}\simeq L$ et 
\begin{eqnarray*}\label{out4}
\mathbb{D}(\mathbb{H})_\Q =\bigoplus_{j\in\Z/n\Z}
\mathbb{D}(\mathbb{H})_{\Q,j}  \underset{\sum_j \Pi^{-j}}{\xrig{\;\;\sim\;\;}}
\mathbb{D} (\mathbb{H})_{\Q,0}^n \simeq L^n
\end{eqnarray*}
Le choix de l'isomorphisme (\ref{out1}) et de ceux qui s'en suivent
n'est pas vraiment nécessaire à la
démonstration mais permet d'identifier les espaces de périodes du coté
Lubin-Tate et Drinfeld à des sous-espaces de $\mathbb{P}^n$.

Rappelons que l'on a des extensions de corps valués  
$$
K| \breve{F} | F |\Qp
$$
et un isomorphisme $L\simeq \breve{F}$.

\subsection{L'application $\M^\D_\infty (K)\ldrt \M^{\LT}_\infty (K)$} \label{dltisog}

Soit $(G,\rho_G,\eta_G)\in \M_\infty^\D (K)$. 
Rappelons (cf. théorème \ref{descDri}) qu'on lui associe un couple
$(\Fil_G,\zeta_G)$ où
$$
\Fil_G\in \Omega (K)
$$
et
$$
\zeta_G\in \Hom_\O (F/\O_F,\mathbb{G}\otimes \O_K/p\O_K)\unp \setminus
\{ 0 \}
$$
Soit 
$$
\zeta_{G*} : K\ldrt \mathbb{D} (\mathbb{G})_{\Q}\otimes K
$$
l'application induite au niveau de l'évaluation des cristaux de Messing
sur l'épaississement $\O_K\twoheadrightarrow \O_K/p\O_K$. Considérons
le composé
$$
\xymatrix@R=5mm{
K\ar[r]^(.3){\zeta_{G*}} & \mathbb{D}(\mathbb{G})_\Q\otimes K
\ar[r]^{\Delta_*^{-1}} & \mathbb{D} (\mathbb{H})_\Q^n\otimes K \ar@{=}[r] &
\bigoplus_{j\in \Z/n\Z} \mathbb{D}(\mathbb{H})_{\Q,j}^n\otimes K \\
1 \ar@{|->}[rrr] &&& (a_{ij})_{1\leq i\leq n,j\in \Z/n\Z}
}
$$
et posons $x_{ij}=\Pi^{-j}a_{ij}\in \mathbb{D}(\mathbb{H})_{\Q,0}\otimes
K\simeq K$. Soit 
$$
X=(x_{ij})_{i,j}\in\text{M}_n (K)
$$

\begin{lemm}\label{enkgf}
L'application $K$-linéaire  de $K^n$ dans lui même
induite par $X$  a
pour image $\Fil_G$.
\end{lemm}
\dem
D'après le théorème \ref{decHTO} et la remarque \ref{HTet} le
sous-module
engendré par l'image de $\zeta_{G*}$ dans
$\mathbb{D}(\mathbb{G})_\Q\otimes K$ est $\bigoplus_{j\in \Z/n\Z}
\Pi^j \Fil_G$. Il en résulte aussitôt que $$K(x_{i 0})_i+\dots + K
(x_{in})_i =Fil_G \subset K^n=\mathbb{D}(\mathbb{G})_{\Q,0}\otimes
K$$
\qed

\begin{defi}
On note $\Fil_H\in \mathbb{P}^n (K)$ l'image de $\,^t X$ dans $K^n$
i.e. le sous-espace engendré par les lignes de $X$.
\end{defi}

Via l'identification $\mathbb{D}(\mathbb{H})_\Q\otimes K \simeq K^n$
si 
$$
\xymatrix@R=5mm{
K \ar[r]^(.3){\zeta_{G*}} & \mathbb{D}(\mathbb{G})_\Q\otimes K
\ar[r]^{\sim} & \mathbb{D}(\mathbb{H})_\Q^n\otimes K \\
1 \ar@{|->}[rr] && (\a_i)_i
}
$$
alors $\Fil_H \in \mathbb{P} ( \mathbb{D} (\Hb))(K)$ est l'image de l'application $K^n\ldrt
\mathbb{D}(\mathbb{H})_\Q\otimes K$ définie par $(\a_i)_i$.
\\

La filtration $\Fil_H$ définit 
 donc d'après la proposition \ref{four_D} un élément 
$\M^{\LT} (K)/\sim$ dont il reste à définir la rigidification du
 module de Tate (cf. théorème \ref{derhj}). 
Considérons le composé
$$
F/\O_F \xrig{\;\; \zeta_G \;\;} \mathbb{G}\otimes \O_K/p\O_K
\xrig{\;\;\Delta^{-1}\;\;} \mathbb{H}^n\otimes \O_K/p\O_K
$$
qui fournit des éléments 
$$
(\zeta_{H,i})_{1\leq i \leq n} \in \Hom_\O (F/\O_F, \mathbb{H}\otimes \O_K/p\O_K)\unp
$$
Pour un entier $i\in \{1,\dots,n\}$ le morphisme induit entre cristaux
de Messing évalués sur $\O_K\twoheadrightarrow \O_K/p\O_K$ est 
$$
\xymatrix@R=5mm{
(\zeta_{H,i})_* : & K \ar[r]^(.3){\zeta_{G*}} & \mathbb{D}
  (\mathbb{G})_\Q\otimes K \ar[r]^{\Delta_*^{-1}} & \mathbb{D}
  (\mathbb{H})_\Q^n\otimes K \ar[r]^{\;\text{proj}_i\;} &
  \mathbb{D}(\mathbb{H})_\Q\otimes K \\
& 1 \ar@{|->}[rrr] &&& \a_i
}
$$
et donc $\forall i\; \text{Im}( (\zeta_{H,i})_*)\subset \Fil_H$. 
\\

Il reste à voir que les $(\zeta_{H,i})$ sont $F$-linéairement indépendants.
Mais si $(\l_i)_i\in F^n$ est tel que $\sum_i \l_i \zeta_{H,i}=0$
alors $(\l_i)_i$ définit une forme linéaire sur $F^n$.
 Et l'égalité $\sum_i \l_i \zeta_{H,i}$ implique sur
l'évaluation des cristaux l'égalité $\sum_i \l_i (a_{ij})_j = 0$ et donc 
$\sum_i \l_i (x_{ij})_j = 0$ ce qui implique d'après le lemme \ref{enkgf}
que la forme linéaire associée à $(\l_i)_i$ s'annule sur $\Fil_G$. Si $(\l_i)_i$ est non-nul alors  $\Fil_G$ est égal au noyau de la forme linéaire associée, est donc défini sur $F$
et contient donc à fortiori une droite $F$-rationnelle.
Donc, puisque $\Fil_G\in \Omega (K)$, $\forall i\; \l_i= 0$.
\\

D'après le théorème \ref{derhj} on en déduit un triplet
$(H,\rho_H,\eta_H)\in \M_\infty^{\LT} (K)$.

\subsection{L'application $\M_\infty^{\LT}(K)\ldrt \M_\infty^{\D} (K)$}

Soit $(H,\rho_H,\eta_H)\in\M_\infty^{\LT} (K)$. Rappelons (cf. théorème
\ref{derhj}) qu'on lui associe un couple $(\Fil_H,\zeta_H)$ où 
$$
\Fil_H\in \mathbb{P}^n (K)
$$
et
$$
\zeta_H=(\zeta_{H,i})_i : F^n\hookrightarrow
\Hom_\O (F/\O_F,\mathbb{H}\otimes \O_K/p\O_K)\unp
$$
Soit 
$$
\xymatrix@R=4mm@C=12mm{
K \ar[rr]^(.37){(\zeta_{H,1*},\dots,\zeta_{H,n*})} &&
\mathbb{D}(\mathbb{H})_\Q^n\otimes K \ar@{=}[r] & \bigoplus_{j\in
  \Z/n\Z} \mathbb{D} (\mathbb{H})_{\Q,j}^n\otimes K \\
1 \ar@{|->}[rrr] &&& (a_{ij})_{i,j}
}
$$
l'application induite sur l´évaluation des cristaux sur
l'épaississement $\O_K\twoheadrightarrow \O_K/p\O_K$. Notons 
$x_{ij}=\Pi^{-j} a_{ij}\in\mathbb{D}(\mathbb{H})_{\Q,0}\otimes K
\simeq K$. Alors,
$$
X=(x_{ij})_{i,j}\in \text{M}_n (K)
$$
\begin{lemm}
L'endomorphisme $K$-linéaire de $K^n$ induit par $\,^tX$ a pour
image $\Fil_H$. 
\end{lemm}
\dem
C'est une conséquence du théorème \ref{decHTO}.
\qed

\begin{defi}
On note $\Fil_G=\text{Im} (X)\in \mathbb{P}^n (K) =
\mathbb{P}(\mathbb{D}(\mathbb{G})_{\Q,0})(K)$ via
$\mathbb{D}(\mathbb{H})_{\Q,0}^n\simeq \mathbb{D}(\mathbb{G})_{\Q,0}$
induit par $\Delta$.
\end{defi}

\begin{prop}
$\Fil_G\in \Omega (K) \subset \mathbb{P}\left (
  \mathbb{D}(\mathbb{G})_{\Q,0}^{V=\Pi}\right )(K)$
\end{prop}
\dem
On a 
$$
\Fil_G = K(x_{i1})_i+\dots + K ( x_{in})_i
$$
Soit $(\mu_j)_j\in K^n$ tel que 
$$
\sum_j \mu_j (x_{ij})_i\in F^n 
$$
Soit alors la forme $F$-linéaire à valeurs dans $F$ définie sur
$V_p(H)$ par  
$$
\ph : \zeta_{H,i}\longmapsto \sum_j \mu_j x_{ij}\in F
$$
Elle définit un élément de $V_p (H)^*$. D'après le corollaire
\ref{DansOmega} dans la suite de Hodge-Tate de $H^D$ 
$$
0\drt \omega_{H^D}^*\otimes K \drt V_p (H)^*\otimes_{F} K \drt
\omega_H\otimes K\drt 0
$$
l'application $V_p(H)^*\ldrt \omega_H\otimes K$ est injective. Or
$$
\left [ V_p (H)^*\otimes_F K\twoheadrightarrow \omega_H\otimes K\right ]
=\left [ \ker \a \hookrightarrow V_p(H)\otimes_F K \right ]^*
$$
où $\a : V_p (H) \otimes_F K \ldrt \mathbb{D}(\mathbb{H})_\Q\otimes K$ est
l'application de matrice $\,^t(a_{ij})_{i,j}$. Mais $\ph_{|\ker \a}=0$ par
définition. Donc $\ph=0$ et $\sum_j \mu_j (x_{ij})_i=0$.
\qed

On obtient donc d'après la proposition \ref{DrP} un couple
$(G,\rho_G)\in \M^{\D}(K)/\sim$. Reste à définir une rigidification de
$G$. Considérons le morphisme composé
$$
\zeta_G : F/\O_F \xrig{\;(\zeta_{H,1},\dots,\zeta_{H,n})\;} 
\mathbb{H}^n\otimes \O_K/p\O_K \xrig{\;\Delta\;}\mathbb{G}\otimes\O_K/p\O_K
$$
On vérifie aussitôt que par définition de $\Fil_G$ le morphisme induit
au niveau de l'évaluation des cristaux a son image contenue dans
$\Fil_G$. D'après le théorème \ref{descDri} on obtient donc un
élément 
$(G,\rho_G,\eta_G)$.
\qed

\subsection{Les deux applications sont inverses l'une de l'autre}

Cela ne pose pas de problème.

\subsection{Retraçage des actions}

\begin{prop}
Dans la bijection entre $\M^{\LT}_\infty (K)$ et $\M_\infty^\D(K)$
si 
$$
(H,\rho_H,\eta_H) \longmapsto (G,\rho_G,\eta_G)
$$
et  $(g,d)\in \GL_n (F)\times D^\times$  alors
$$
(g,d).(H,\rho_H,\eta_H) \longmapsto (\,^t g,d^{-1}) .(G,\rho_G,\eta_G)
$$
\end{prop}
\dem
Cela ne pose aucun problème.
\qed

\subsection{Bijection entre les points de l'espace de Berkovich
  associé}

Soit $K\hookrightarrow K'$ une extension isométrique de corps valués
complets pour une valuation de rang $1$. On vérifie alors aussitôt que
le diagramme suivant est commutatif
$$
\xymatrix{
\M_\infty^{\LT} (K) \ar[d] \ar[r]^\sim & \M_\infty^{\D} (K) \ar[d] \\
\M_\infty^{\LT} (K') \ar[r]^\sim & \M_\infty^{\D} (K')
}
$$
Pour $*\in \{ \LT,\D\}$ posons 
$$
\left |\M^*_\infty \right |= \coprod_K \M_\infty^* (K) /\sim
$$
où pour $x\in \M^*_\infty (K_1)$ et $y\in \M^*_\infty (K_2)$  
$\; x\sim y$ ssi il existe une extension valuée comme précédemment
$$
\xymatrix@R=2mm{ K_1 \ar@{^(->}[rd] \\
& K_3 \\
K_2 \ar@{^(->}[ru]
}
$$
telle que $x$ et $y$ aient même image dans $\M^*_\infty (K_3)$.

Cette classe d'équivalence est bien définie au sens où c'est un ensemble. Cela
résulte de l'existence de modèles entiers de nos espaces qui implique
que l'on peut se limiter à des corps $K$ de cardinalité bornée. On a donc
une bijection $\GL_n (F)\times D^\times$-équivariante 
$$
\left | \M^{\LT}_\infty \right | \iso \left | \M^{\D}_\infty  \right |
$$
Le choix de modèles entiers de nos espaces permet de munir ces
ensembles d'une structure d'espace topologique localement compact et
l'existence de l'isomorphisme au niveau de ces modèles entiers
(\cite{iso4}) impliquera que la bijection précédente est un
homéomorphisme. 

\section{La matrice  $X$}\label{kfyygpzr5rrppn}

Les triplets $(H,\rho_H,\eta_H)$ et $(G,\rho_G,\eta_G)$ 
qui se correspondent ont en commun la
matrice de rang $n-1$  $\;X\in\text{M}_n (K)$
où comme précédemment on identifie
$$
\DH_\Q^n\otimes K \simeq \bigoplus_{j\in \Z/n\Z} \DH_{\Q,j}^n\otimes K
\underset{\sim}{\xrig{\oplus_j \Pi^{-j}}} \bigoplus_{j\in \Z/n\Z}
\DH_{\Q,0}^n\otimes K \simeq M_n (K)
$$
 Celle-ci vérifie
\begin{itemize}
\item Les lignes de $X$ engendrent $\Fil_H\in \mathbb{P}^n (K)$
\item Les colonnes de $X$ engendrent $\Fil_G\in\Omega (K)$
\end{itemize}

Supposons maintenant pour simplifier que $F=\Qp$. Alors, $X$ possède
un relèvement 
$$
X_{cris}\in \text{M}_n ( B^+_{cris}(\O_K))
$$
i.e. via $\theta : B^+_{cris} (\O_K)\twoheadrightarrow
\O_{\widehat{K}}$ on a $\theta (X_{cris})=X$.
 Ce relèvement est défini
de manière similaire à $X$ en évaluant les cristaux sur
$A_{cris}(\O_K)\twoheadrightarrow \O_{\widehat{\overline{K}}}/p  \O_{\widehat{\overline{K}}}$ au
lieu de $\O_K\twoheadrightarrow \O_K/p\O_K$ (cf. les propositions  \ref{docyum1} et
\ref{docyum2} de la section \ref{ruhgp}
pour plus de détails).

De plus,
$$
\det (X_{cris}) \in \Qp^\times .t 
$$
où la valuation de l'élément de $\Qp^\times$ est liée aux hauteurs de
$\rho_H,\rho_G$ et $\Delta$ (cf. la section suivante). 
Si 
$$\Phi=\left (
\begin{matrix}
0 & p &  & \\
0 & 0& \ddots &  \\
\vdots  &  &\ddots & p \\
 1 &0 & & 0  
\end{matrix}\right )
$$
est la matrice du Frobenius de $\mathbb{D}(\mathbb{H})_\Q$
relativement à la base définie par l'isomorphisme 
$\mathbb{D}(\mathbb{H})_{\Q,0}\simeq L$ et $\DH_\Q\simeq \oplus_j
\DH_{\Q,j} \underset{\sim}{\xrig{\sum_j \Pi^{-j}}} \oplus_j \DH_{\Q,0}$
\begin{eqnarray}\label{eqslf}
\boxed{
\ph (X_{cris}) \,^{t}\Phi = p X_{cris} }
\end{eqnarray}
où $\ph (X_{cris})$ est obtenue en appliquant le Frobenius cristallin
à tous les éléments de la matrice $X_{cris}$. 
On peut montrer la proposition suivante :
\begin{prop}
Les espaces $\M_\infty^{\LT} (K)$ et $\M_\infty^{\D}(K)$ s'envoient
surjectivement sur l'ensemble des matrices $X\in\text{M}_n (K)$ de rang
$n-1$ possédant un relèvement  $X_{cris}$ à $B^+_{cris}$ de
déterminant non-nul et vérifiant  
l'équation (\ref{eqslf}). De plus les fibres de cette application en $X$
sont
en bijection avec l'ensemble des relèvements de $X$
de déterminant non-nul satisfaisant à l'équation
(\ref{eqslf}).

 Cette ensemble de matrices est muni d'une action de
$\GL_n (F)\times D^\times$ de la façon suivante : l'isomorphisme
 $\mathbb{D}(\mathbb{H})_\Q\simeq L^n$ induit un plongement $D^\times
 \hookrightarrow \GL_n (L)$. Alors 
$$\forall (g,d)\in \GL_n (F)\times
 D^\times \;\;\forall X \;\;\;\;\;\; (g,d).X= \,^tg X d^{-1}  
$$ 
Les applications précédentes sont compatibles à cette action.
 \end{prop}

$$
\xymatrix@C=6mm@R=22mm{
\M_\infty^{\D} (K) \ar[rrd]\ar[dd]_{Periodes} \ar[rrrr]^{\sim} &&&&
\M_\infty^{\LT}(K) \ar[lld]  \ar[dd]^{Periodes} \\
 && \{X\in \text{M}_n (K)\; |\; \text{rg} X =n-1\; \}
\ar[lld]^{X\mapsto \text{Im} \,^t X} \ar[rd]^{X\mapsto \text{Im} X} \\
\mathbb{P}^n (K) && &  \mathbb{P}^n (K) & \Omega(K) \ar@^{_(->}[l]
}
$$
\dem
Il s'agit d'une application du théorème \ref{deall} énoncé dans
l'appendice.
\qed

\begin{rema}
Dans la proposition précédente on peut enlever  ``de déterminant
non-nul'' si l'on se restreint aux $X$ telles que $\text{Im}\, X \in \Omega (K)$.
\end{rema}

\begin{rema}
Pour les corps locaux de caractéristique positive
A.Genestier et V.Lafforgue donnent une interprétation géométrique comme points d'un certain schéma formel  
de l'analogue en caractéristique positive de l'ensemble des matrices $X_{cris}$ solutions de l'équation précédente.  Il n'existe pas de telle interprétation en caractéristique zéro.
\end{rema}

\section{L'isomorphisme conserve le degré }\label{ruhgp}

Il s'agit de démontrer que pour un choix convenable de $\Delta$ le diagramme suivant commute 
$$
\xymatrix{
\M_\infty^\D (K) \ar[rd] \ar[rr]^\sim &&  \M_\infty^{\LT} (K) \ar[ld] \\
& \Z
}
$$
où les applications vers $\Z$ sont les applications hauteur renormalisées  
 des quasi-isogénies universelles. 
Cela implique en particulier que les tours de Lubin-Tate et Drinfeld
classiques, les fibres en hauteur zéro, sont isomorphes. 

Plus précisément, si $f$ est une quasi-isogénie entre $\O$-modules
formels on note $\text{ht}_\O (h) = \text{ht} (f) /[F:\Qp]$. Alors
l'application $\M_\infty^\D (K) \ldrt \Z$ est 
$$
(G,\rho_G,\eta_G)\longmapsto \text{ht}_\O (\rho_G) /n
$$
et l'application $\M_\infty^{\LT} (K)\ldrt \Z$ est 
$$
(H,\rho_H,\eta_H)\longmapsto \text{ht}_\O (\rho_H)
$$

On va appliquer le théorème \ref{det_periodes_unite} de l'appendice.
Afin de simplifier on se limite au cas $F=\Qp$, le cas général étant
laissé au lecteur.
\\

Supposons que l'isomorphisme  $\DH_{\Q,0} \simeq L$ fixé au début de  la
section \ref{coeulop} provienne d'un
isomorphisme $\DH_0\simeq \O_L$. 

\begin{prop}\label{docyum1}
 Soit $(H,\rho_H,\eta_H)\in \M_\infty^{\LT} (K)$.
Soit $\xi : \Qp/\Zp \ldrt \Hb^n\otimes \O_K/p\O_K$ associé à $\eta_H$, $\rho_H$ et la base canonique
de $\Zp^n$. Soit 
$$
\xi_* : B^+_{cris} \ldrt \DH_\Q^n\otimes_L B^+_{cris}
$$
le morphisme associé sur l'évaluation des cristaux. 
 Soit $X_{cris} \in
 \M_n \left ( B^+_{cris} (\O_K)\right ) $ la matrice associée via
$$
\DH_\Q^n\otimes_L B^+_{cris} = \bigoplus_{j\in\Z/n\Z} \DH_{\Q,j}^n
\otimes_L B^+_{cris}
 \underset{\sim}{\xrig{\; \oplus_j
    \Pi^{-j}\; }} \bigoplus_{j\in \Z/n\Z} \DH_{\Q,0}^n\otimes_L B^+_{cris}
$$
Alors $\det \left ( X_{cris}\right )  = \l .t \in \Qp^\times .t$ et de plus
$$
v_p (\l) = -\text{ht} (\rho_H) - \frac{1}{2} n (n-1)
$$ 
\end{prop}
\dem
Soit $M_0 \subset \DH_\Q$ le module de Dieudonné de $H\otimes_{\O_K} k$
$$
M_0 = \mathbb{D} (\rho_H^{-1}) \left ( \mathbb{D} (H_k)\right )
$$
D'après le théorème \ref{det_periodes_unite} de l'appendice, dans
l'isomorphisme
$$
T_p (H) \otimes_{\Zp} B^+_{cris} \simeq M_0 \otimes_{\O_L} B^+_{cris}
$$
on a 
$$
\bigwedge^n T_p (H)\otimes \O_L = t. \bigwedge^n M_0
$$
De plus, $[M_0 : \DH ] = \text{ht} (\rho_H)$. Via l'isomorphisme 
$\DH_\Q \simeq \oplus_{j\in \Z/n\Z} \DH_{\Q,0}$ le réseau $\DH$
s'envoit sur 
$$
\bigoplus_{j\in\Z/n\Z} \Pi^{-j} \DH_0
$$
Soit $W$ l'image de $T_p (H)$ dans $\oplus_j \DH_0 \otimes_{\O_L}
B^+_{cris}$. Alors
$$
\bigwedge^n W\otimes \O_L = p^{-\text{ht} (\rho_H)-\frac{1}{2} n (n-1)}
t. \DH_0$$
D'où le résultat.
\qed

\begin{prop}\label{docyum2}
Soit $(G,\rho_G,\eta_G)\in\M_\infty^\D (K)$. Soit $\zeta=
\rho_G^{-1}\circ \eta_G (1) : \Qp/\Zp \ldrt \Gb\otimes \O_K/p\O_K$ et
$\xi =\Delta^{-1} \circ \zeta : \Qp/\Zp \ldrt \Hb^n\otimes
\O_K/p\O_K$. Soit $X_{cris}$ la matrice associée comme dans la
proposition précédente. Alors $\det ( X_{cris})= \l . t \in \Qp^\times
t$ et 
$$
v_p (\l) = -\text{ht} (\rho_G)/n - \text{ht} (\Delta)
$$
\end{prop}
\dem
La démonstration est identique à celle de la proposition
précédente. On applique le théorème \ref{det_periodes_unite} de
l'appendice. La matrice des périodes de $G$ est $n$-copies de la
matrice $X_{cris}$. On en déduit que $\det (X_{cris})^n$ est donné par
la valeur annoncée.
\qed

\begin{coro}
Si $\text{ht}_\O (\Delta) = \frac{1}{2} n (n-1)$ alors la bijection entre
$\M_\infty^{\LT} (K)$ et $\M_\infty^\D (K)$ respecte les hauteurs normalisées.
\end{coro}

\begin{coro}
Soit $(H,\rho_H,\eta_H)\in \M_\infty (K)$ tel que $\text{ht}
(\rho_H)=0$. Alors, le couple $(G,\rho_G)$ associé dans l'espace de
Rapoport-Zink sas  niveau est  l'élément de $\Omega (K)
\subset \M^{\D} (K)$ donné 
par le dual de la décomposition de Hodge-Tate de $H$ couplé à $\eta_H$.
Plus généralement, soit $(G,\rho_G)$ le point de $\Omega (K)$ associé à la décomposition de  Hodge-Tate de $H$ et à $\eta_H$. Le point de $\M^{\D} (K)$ 
 associé par l'isomorphisme
de Faltings est alors 
$$
(G,\rho_G\circ \Pi^{\text{ht} (\rho_H)})
$$
\end{coro}

\section{Un point de vue différent sur la bijection} \label{ptovue}

Nous allons redémontrer la bijection précédente d'un point de vue ``dual''.
 En particulier, au lieu de considérer des filtrations $\Fil_H\in
 \mathbb{P}^n, \Fil_G\in \Omega$ nous considérerons plutôt les
 quotients $K^n\twoheadrightarrow K^n/\Fil_H, K^n\twoheadrightarrow
 K^n/\Fil_G$. Ce point de vue se prête mieux lorsque l'on travaillera
 sur une base quelconque (c'est à dire plus sur un point comme dans
 cet article)  puisque
 que pour un fibré vectoriel $\mathcal{E}$ il est plus
 commode de définir $\mathbb{P} (\mathcal{E})$ comme classifiant les
 quotients localement libres de rang $1$
 $\mathcal{E}\twoheadrightarrow \mathcal{L}$ (en tous cas cela est
 plus facile à platifier par éclatements). En effet, dans le cas
 d'une base quelconque l'approche des sections précédentes nous conduirait à définir $\Fil_H$
 comme sous-module engendré par certaines sections ce qui est moins
 commode. De plus la matrice définissant ces sous-modules serait une
 section de  $\DH_\Q^n\otimes
 \O_\X \unp\simeq M_n  ( \O_\X \unp  )$ où $\X$ est  l'espace de 
 Lubin-Tate ou de Drinfeld en niveau infini. Il faudrait donc procéder
 à de nouveaux éclatements afin de rendre l'espace engendré par les
 lignes et celui par les colonnes localement facteur direct entier 
(cette dernière justification est quelque peu hypocrite puisque de
toutes façons dans la démonstration finale on devra à un endroit rendre
entier l'application des périodes). Le point de vue qui suit permet de
construire directement une application entière du sommet de la tour de
Lubin-Tate (resp. la tour de Drinfeld) vers le schéma formel de Drinfeld
(resp. $\widehat{\P}^{n-1}$).

Commençons par faire le lien entre le point de vue précédent et celui
qui va suivre.

\subsection{Identification de $K^n\twoheadrightarrow K^n/\Fil_H$ avec
  l'application de Hodge-Tate de $G^D$}

Soit $(G,\rho_G,\eta_G)\in \M^{\D}_\infty (K)$. Rappelons qu'on a
défini une matrice $X$ telle que $\text{Im} X =\Fil_G$ et
$\text{Im}\,^t X =\Fil_H$. Tout repose sur une bidualisation,
l'exactitude de la suite de Hodge-Tate ainsi que
sur un analogue de
l'énoncé d'algèbre linéaire suivant : soit $u:E_1\drt E_2$ une
application linéaire entre deux $K$-espaces vectoriels de dimension
finie, il y a alors une identification canonique 
$$
\left [ E_1^*\twoheadrightarrow E_1^*/\text{Im} \,^t u\right ]
\simeq \left [ \ker u \hookrightarrow E_1 \right ]^*
$$
L'analogue est le suivant :  pour $E$ un $D\otimes_F K$-module de type
fini (à gauche ou à droite) posons 
$$
\Gamma (E) =\Hom_{D\otimes_F L}(E, \mathbb{D}(\mathbb{H})_{\Q})
$$
un $K$-espace vectoriel de dimension finie. Alors, 
\begin{eqnarray}\label{ergjk}
\forall u : E_1\drt E_2\;\;\; \left [
\Gamma (E_1)\twoheadrightarrow \Gamma (E_1)/ \text{Im} \,^\Gamma u
\right ] \simeq \Gamma \left [\ker u \hookrightarrow E_1\right ]  
\end{eqnarray}
Notons maintenant pour $W$ un $K$-e.v. $\Phi (W) =W^*$ le dual
usuel. Si $E$ est un $D\otimes_L K$-module de type fini à gauche,
resp. à droite, alors $\Phi (E)$ est naturellement un $D\otimes_L
K$-module de type fini à droite, resp. à gauche.

 Il y a alors un isomorphisme naturel de bidualité pour $E$ comme ci-dessus 
$$
 \GG\circ \Phi (E) \simeq \mathbb{D} (\mathbb{H})_{\Q,0} \otimes_L E_0
$$
où  comme d'habitude $E_0 = \{x\in E\; |\; \forall a\in F_n\;  a\otimes
1. x = 1\otimes a.x \;\}$. 

\begin{lemm}
Soit l'application $D\otimes_F K$-linéaire composée
$$
\xymatrix@C=14mm{
u : D\otimes_F K  \ar[r]^{\eta_G}_{\sim} & V_p (G)\otimes_F K
\ar[r]^{\a_G} &
\mathbb{D}(\mathbb{G})_\Q\otimes_L K
}
$$
Procédons aux identifications suivantes :
\begin{eqnarray*}
\Gamma (D\otimes_F K) & \iso & \mathbb{D}(\mathbb{H})_\Q\otimes_L K \\
f &\longmapsto & f (1\otimes 1)
\end{eqnarray*}
et
\begin{eqnarray*}
\Gamma (\mathbb{D} (\mathbb{G})_\Q\otimes_L K) &\iso & \Hom_L (
\mathbb{D}(\mathbb{G})_{\Q,0},\mathbb{D}(\mathbb{H})_{\Q,0})\otimes_L
K \simeq K^n \\
f &\longmapsto & f_{|\mathbb{D}(\mathbb{G})_{\Q,0}\otimes K}
\end{eqnarray*}
via l'isomorphisme (\ref{out1}) et $\Delta$. Alors 
$$
\,^\Gamma u : K^n\ldrt \mathbb{D} ( \mathbb{H})_\Q\otimes K \;\;\;\;(\simeq K^n)
$$
s'identifie à $\,^t X$ et donc $\text{Im} (\,^\Gamma u) =\Fil_H$. 
\end{lemm}
\dem
Il suffit de retracer les différentes identifications.
\qed

\begin{rema}
\'Etant donné que le module de Tate de $G$ est trivialisé et que
$\det V_p (G)\simeq \Qp (1)$ on en déduit que
$\Qp (1)$ est trivialisé et que donc l'application de
Hodge-Tate de $G^D$ est définie sur $K$ : $\a_{G^D} : V_p (G^D)
\ldrt \omega_G\otimes K$. En d'autres termes $\Qp
(\mu_{p^{\infty}})\subset K$.  
\end{rema}

\begin{coro}
Il y a des identifications canoniques 
\begin{eqnarray*}
\left [ \mathbb{D}(\mathbb{H})_\Q\otimes K \twoheadrightarrow
  \mathbb{D}(\mathbb{H})_\Q\otimes K/\Fil_H \right ] &\simeq &
\left [ (V_p (G^D)\otimes_F K)_0  \overset{\a_{G^D}}{\twoheadrightarrow}
  \omega_{G,0}\otimes K  \right ] (-1)\otimes_L \mathbb{D} (\mathbb{H})_{\Q,0} \\
&\simeq & \left [\mathbb{D}(\mathbb{H})_\Q\otimes K\twoheadrightarrow 
\omega_{G,0}(-1)\otimes_L  \mathbb{D}(\mathbb{H})_{\Q,0}\right ]
\end{eqnarray*}
où $\a_{G^D} $ est l'application de Hodge-Tate de $G^D$.
\end{coro}
\dem
Soit 
$$
\left [ E_1 \xrig{\;\a_G\;} E_2 \right ] = \left [ V_p (G)\otimes_F K
  \xrig{\;\a_G\;} \mathbb{D} (\mathbb{G})_\Q\otimes_L K \right ]
$$
D'après le lemme précédent
$$
\left [ \mathbb{D} (\mathbb{H})_\Q\otimes K \twoheadrightarrow
  \mathbb{D} (\mathbb{H})_\Q\otimes K/\Fil_H \right ]
\simeq 
\left [ \GG (E_1) \overset{\,^\GG \a_G}{\twoheadrightarrow} \GG
  (E_1)/\text{Im} (\,^\GG \a_G ) \right ]
$$
qui d'après l'identification (\ref{ergjk}) est isomorphe à
$$
\GG \left [ \ker \a_G \hookrightarrow E_1 \right ]
$$
Mais d'après la proposition \ref{relhtO} (décomposition de Hodge-Tate)
la suite
$$
0 \ldrt \omega_G^*\otimes K(-1) \xrig{\;\,^\Phi\a_{G^D} (-1)\;} V_p
(G)\otimes_F K \xrig{\; \a_G\; } \mathbb{D} (\mathbb{G} )_\Q\otimes_L
K 
$$
est exacte. Donc,
\begin{eqnarray*}
\simeq \GG \left [ \ker \a_G \hookrightarrow E_1 \right ]& \simeq&
\GG\circ \Phi \left [ V_p (G^D)\otimes_F K \twoheadrightarrow
  \omega_G\otimes K \right ](-1) \\
&\simeq & \left[ \left ( V_p (G^D)\otimes_F K \right )_0 (-1) \twoheadrightarrow
  \omega_{G,0}\otimes K(-1) \right ]\otimes_L\mathbb{D} (\mathbb{H})_{\Q,0}
\end{eqnarray*}
Le dernier isomorphisme dans l'énoncé du corollaire provient de la
rigidification
$\eta : D\iso V_p (G)$ qui induit 
$$
\GG\circ \Phi \left ( V_p (G)\otimes_F K\right ) \iso \mathbb{D}
  (\mathbb{H})_{\Q}\otimes K
$$
\qed

\subsection{Identification de $K^n\twoheadrightarrow K^n/\Fil_G$ avec
  l'application de Hodge-Tate de $H^D$}

Soit $(H,\rho_H,\eta_H)\in \M_\infty^{\LT} (K)$. On procède comme dans
la section précédente. Pour $E$ un $K$-espace vectoriel posons
$$
\Psi (E) = \Hom_L (E,\DH_{\Q,0} )
$$
\begin{lemm}
Soit l'application composée 
$$
 \xymatrix{
v:K^n \ar[r]^(.4){\eta_H}_(.4){\sim} & V_p (H)\otimes_F K \ar[r] &
 \DH_\Q\otimes_L K
}
$$
Procédons aux identifications suivantes :
\begin{eqnarray*}
\Psi (K^n)&\simeq& \DH_{\Q,0}^n\otimes K \simeq \DG_{\Q,0}\otimes
K\text{   via } \Delta \\
\Psi ( \DH_\Q\otimes_L K) &\simeq & \Hom ( \DH_\Q,\DH_{\Q,0})\otimes K
\simeq K^n
\end{eqnarray*}
où la second identification utilise
$$
\DH_\Q = \bigoplus_{j\in\Z/n\Z} \DH_{\Q,j} \underset{\sim}{\xrig{\oplus_i
    \Pi^{-i}}}  \DH_{\Q,0}^n\simeq L^n
$$
Alors $\,^\Psi v : K^n \ldrt \DG_{\Q,0} \otimes K$ s'identifie à la
matrice
 $X$ et donc $\text{Im} (\,^\Psi v)=\Fil_G$. 
\end{lemm}

\begin{coro}
Il y a des identifications 
\begin{eqnarray*}
\left [ \DG_{\Q,0}\otimes K \twoheadrightarrow  \DG_{\Q,0}\otimes
  K / \Fil_G \right ] & \simeq &  \left [ V_p (H^D)
  \xrig{\;\a_{H^D}\;}\omega_H\otimes K\right ](-1) \otimes_L \DH_{\Q,0} \\
 & \simeq & \left [ \DG_{\Q,0}\otimes K \twoheadrightarrow
   \omega_H\otimes K(-1) \otimes_L \DH_{\Q,0} \right ]
\end{eqnarray*}
\end{coro}
\dem
Soit 
$$
\left [ E_1\xrig{\;\a_{H}\;} E_2\right ] = \left [ V_p (H)\otimes_F K
  \xrig{\;\a_H\;} \DH\otimes_L K \right ]
$$
Alors
\begin{eqnarray*}
\left [ \DG_{\Q,0}\otimes K \twoheadrightarrow \DG_{\Q,0}\otimes
  K / \Fil_G\right ] & \simeq & \left [ \Psi (E_1) \twoheadrightarrow
  \Psi ( E_2) \right ] \\ 
& \simeq & \Psi \left [ \ker \a_H \hookrightarrow E_1 \right ] \\
& \simeq & \Psi\circ \Phi \left [ V_p (H^D)\otimes_F K(-1)
  \xrig{\;\a_{H^D} (-1)\;} \omega_H\otimes K(-1) \right ] \\
& \simeq & \left [ V_p (H^D)\otimes_F K
  \overset{\a_{H^D}}{\twoheadrightarrow} \omega_H\otimes K \right
](-1) \otimes_L \DH_{\Q,0} 
\end{eqnarray*}
où on a utilisé l'exactitude de la suite de Hodge-Tate pour $H$
(proposition \ref{relhtO}).
\qed

\subsection{L'application $\M^\D_\infty (K)\ldrt \M^{\LT}_\infty (K)$} \label{dltisog2}

Soit $(G,\rho_G,\eta_G)\in \M_\infty^\D (K)$. 

\subsubsection{Première étape : on  met une structure d'isocristal sur le
  module de Tate}

Posons  
$$
N=\Hom_D ( V_p (G), \mathbb{D} (\Hb)_\Q ) 
$$
un isocristal relativement à l'extension $L|F$ si on  le munit de $\ph : N\iso N$ défini par $\ph.f=\ph \circ f$.
 
\subsubsection{Deuxième étape : la rigidification du module de Tate
  induit une  rigidification de l'isocristal}
 
L'isomorphisme de $D$-modules $\eta_G : D\iso V_p (G)$ induit un isomorphisme d'isocristaux
$$
\mathbb{D} (\Hb)_\Q\iso ( N,\ph)
$$

\subsubsection{Troisième étape : la filtration de Hodge-Tate induit une
  filtration du module de Dieudonné} 

Considérons l'application de Hodge-Tate de $G^D$ tordue par $K(-1)$
$$
(\mathcal{S}) = \left [ V_p (G)^*\otimes_F K \overset{\a_{G^D}(-1)}{\twoheadrightarrow}
  \omega_G\otimes K(-1) \right ] = \bigoplus_{j\in \Z/n\Z} \left [ 
\left ( V_p (G)^*\otimes_F L\right )_j\otimes_L K \twoheadrightarrow
  \omega_{G,j}\otimes K(-1) \right ]
$$
Tensorisons cette suite
$$
(\mathcal{S})\otimes_{D\otimes_F L} \DH_\Q \simeq \left [ N\otimes_F
  L\twoheadrightarrow \omega_{G,0}(-1)\otimes_L \DH_{\Q,0} \right ]
$$
qui fournit donc via la deuxième étape un élément
$$
\Fil_H\in \mathbb{P} \left ( \DH_\Q \right ) (K)
$$
On obtient ainsi un couple $[(H,\rho_H)]\in \M^{\LT} (K)/\sim$.

\subsubsection{Quatrième étape :   construction d'éléments dans le module de Tate de $H$}

Construisons un morphisme 
$$
\Psi : \Hom_ {\O_D} (\Gb, \Hb)\unp \ldrt V_p (H)
$$
Soit 
$$
\zeta_G \in \Hom (F/\O_F, \Gb\otimes \O_K/p\O_K)\unp
$$
comme dans le théorème \ref{descDri}. Soit $f\in \Hom_{\O_D}
(\Gb,\Hb)\unp$. Considérons $$f\circ \zeta_G \in \Hom (F/\O_F, \Hb
\otimes \O_K/p\O_K)\unp$$
 Pour le morphisme induit sur les cristaux
évalués sur l'épaississement $\O_K\twoheadrightarrow \O_K/p\O_K$
$$
(f\circ \zeta_G)_* : K \ldrt \DH_\Q\otimes K
$$
montrons que $(f\circ \zeta_G)_* (1) \in \Fil_H$. 
\\
Via $\eta_G : D\iso V_p (G)$ le morphisme $f_* : \DG_\Q\otimes K \ldrt
\DH_\Q\otimes K$ s'identifie à
$$
\left ( V_p (G)^*\otimes_F K \right )\otimes_{D\otimes L} \left [ 
\DG_\Q \xrig{\; \mathbb{D} (f)\;}  \DH_\Q \right ]
$$ 
On doit montrer que $\zeta_{G*} (1) \in \left ( V_p (G)^* \otimes
  K\right )\otimes_{D\otimes L} \DG_\Q$ s'envoit sur zéro par la
composée donnée par la ligne  pointillées diagonale dans le diagramme suivant 
$$
\xymatrix@C=4mm@R=6mm{\left ( V_p(G)^*\otimes K \right
  )\ar[dd]^{\a_{G^D}(-1)}  &  \otimes_{D\otimes L} 
    [ \DG_\Q \ar[r]^(0.6){\mathbb{D}(f)} \ar@<3ex>[dd]^{Id} \ar@{-->}@<1.5ex>[rdd] &
    \DH_\Q\ar[dd]^{Id}  ] \\
 & &  \\
\omega_G\otimes K(-1) &     \otimes_{D\otimes L} 
    [ \DG_\Q \ar[r]^{\mathbb{D}(f)} &  \DH_\Q  ]
}
$$
Il résulte de ce diagramme qu'il suffit de vérifier que l'image de
$\zeta_{G*} (1)$  par l'application 
$$
\a_{G^D} (-1) \otimes Id :  \left ( V_p(G)^*\otimes K \right
  )   \otimes_{D\otimes L} \DG_\Q \ldrt \omega_G\otimes K(-1)    \otimes_{D\otimes L} 
     \DG_\Q 
$$
est nulle.
Mais si $\iota : \widetilde{\omega}_{G^D} \otimes K \hookrightarrow \DG_\Q\otimes
K$ et $\,^t \a_{G} : \widetilde{\omega}^*_{G^D}\otimes K \ldrt V_p (G)^*\otimes K$
est l'application de gauche dans la suite de Hodge Tate de $G^D$
tordue par $K(-1)$, $\,^t \a_{G} \in V_p
(G)^*\otimes K \otimes \omega_{G^D} $, alors
$$
\zeta_{G*}(1) = (Id \otimes \iota) \left (\,^t\a_G \right )
$$
Le résultat se déduit donc du fait que dans la suite de Hodge-Tate de
$G^D$ la composée des deux applications est nulle : $\,^t \a_{G} \circ
\a_{G^D} (-1) = 0$. 
On a donc bien défini l'application $\Psi$.

\subsubsection{Sixième étape : l'application $\psi$ est un isomorphisme}

Il suffit de démontrer qu'elle est injective. Soit donc $f$ tel que
$\psi (f)=0$. On vérifie alors sur l'évaluation des cristaux que cela
implique que pour $\mathbb{D} (f) : \DG_\Q\ldrt \DH_\Q$ 
$$
\mathbb{D} (f)_{|\widetilde{\omega}_{G^D}\otimes K} =0
$$
Mais 
$$
\xymatrix@R=2mm{
\mathbb{D} (f) \in \Hom_{D,\ph} \left ( \mathbb{D} (\Gb)_\Q,\mathbb{D}
  (\Hb)_\Q\right ) \ar[r]^\sim & \Hom\left (\mathbb{D}
  (\Gb)_{\Q,0}^{V=\Pi},\DH_{\Q,0}^{V=\Pi}\right ) \\
 h \ar@{|->}[r]  & h_{|\DG_{\Q,0}^{V= \Pi}}
}
$$
or 
$$
\left [  \mathbb{D} (\Gb)_{\Q,0}^{V=\Pi}\otimes K \twoheadrightarrow
  \DG_{\Q,0} \otimes K/ \widetilde{\omega}_{G^D,0}\otimes K
 \right ]\in \Omega (K) \subset \mathbb{P}
\left ( \mathbb{D} (\Gb)_{\Q,0}^{V=\Pi}\right ) (K)
$$
Donc, $\mathbb{D} (f)_{|\omega_{G^D}\otimes K}=0 \impl \mathbb{D} (f)=0$ et donc $f=0$.
\\

La quasi-isogénie $\Delta$ induit un isomorphisme 
$$
F^n \iso \Hom_{\O_D} (\Gb,\Hb)\unp
$$
et
on obtient donc une rigidification $\eta_H$ de $[(H,\rho_H)]$ ce qui
détermine le triplet $(H,\rho_H,\eta_H)$ d'après le théorème \ref{derhj}.

\subsection{L'application $\M_\infty^{\LT} (K) \ldrt \M_\infty^\D (K)$}

Soit maintenant $(H,\rho_H,\eta_H)\in \M_\infty^{\LT} (K)$. 

\subsubsection{Première étape : On met une structure d'isocristal sur le module de Tate}

Soit 
$$
N=\Hom_F ( V_p (H), \mathbb{D} (\Hb)_\Q)
$$
qui est muni d'une structure d'isocristal en posant $\ph.f= \ph\circ f$. Notons que cet isocristal est muni d'une action de $D$
puisque c'est le cas de $\mathbb{D} (\Hb)_\Q$.

\paragraph{Deuxième étape : $\eta_H$ rigidifie l'isocristal}

L'isomorphisme $\eta_H:F^n \iso V_p (H)$ couplé à la quasi-isogénie $\Delta$ 
 induit un isomorphisme d'isocristaux munis d'une action de $D$ 
$$
\mathbb{D} ( \Gb)_\Q \simeq \mathbb{D} (\Hb)_{\Q}^n \iso N
$$

\subsubsection{Troisième étape : la filtration de Hodge-Tate induit une
  du module  de Dieudonné}

Considérons l'application de Hodge-Tate de $H^D$ tordue par $K(-1)$ 
$$
(\mathcal{T})=\left [ V_p (H)^*\otimes K \twoheadrightarrow \omega_H
  \otimes K(-1) \right ]
$$
Après application de $-\otimes_F \DH_{\Q}$ on obtient via la deuxième
partie une filtration $D$-invariante
$$
(\mathcal{T})\otimes_L \DH_\Q \simeq \left [ \DG_{\Q}\otimes K
  \twoheadrightarrow \omega_H\otimes K(-1) \otimes \DH_\Q \right ]
$$
La partie ``indice zéro'' de cette application est obtenue par 
$$
(\mathcal{T})\otimes_L \DH_{\Q,0} \simeq \left [ \DG_{\Q,0}\otimes K
  \twoheadrightarrow \omega_H\otimes K(-1) \otimes \DH_{\Q,0} \right ]
$$
qui définit $\Fil_G\in \mathbb{P} (\DG_{\Q,0}) (K)$. La quasi-isogénie
$\Delta$ couplée à $\eta_H$ induit un isomorphisme 
$$
V_p (H)^*\otimes_F \DH_{\Q,0}^{V^{-1}\Pi} \iso \DG_{\Q,0}^{V^{-1}
  \Pi}
$$
D'après le corollaire \ref{DansOmega} $(\mathcal{T}) \in \Omega (K)
\subset \mathbb{P} (V_p (H)^*)(K)$. Donc 
$$
\Fil_G \in \Omega (K)\subset \mathbb{P} \left ( \DG_{\Q,0}^{V^{-1}\Pi}
\right ) (K)
$$
On obtient donc ainsi un couple $[(G,\rho_G)]\in \M_\infty^{\D} (K)/\sim$.

\subsubsection{ Quatrième étape :  Construction d'éléments dans le module
  de Tate de $G$}

A la rigidification $\eta_H: F^n \iso V_p (H)$ est associée un élément
$$
\xi \in \Hom (F/\O_F,\Hb^n\otimes \O_K/p\O_K)\unp
$$
qui composé avec $\Delta$ fournit $\zeta_G\in \Hom (F/\O_F,\Gb\otimes
\O_K/p\O_K)\unp$.
Il s'agit de montrer que sur l'évaluation des cristaux
$$
\zeta_{G*}(1)   \in \bigoplus_{j\in \Z/n\Z} \Pi^j \Fil_G\subset
\DG_\Q\otimes K
$$
Pour cela il suffit de vérifier que 
$$
\xi_* (1) \in \ker \left [ \DH_\Q^n\otimes K \overset{\eta_H}{\simeq}
  V_p (H)^*\otimes K \otimes \DH_\Q \twoheadrightarrow \omega_H\otimes
  K(-1) \otimes \DH_\Q \right ]
$$
Or cela résulte de ce que $\xi_* (1) \in V_p (H)^*\otimes K\otimes
\DH_\Q$ est donné par 
$$
V_p (H) \xrig{\; \a_H \;} \DH_\Q\otimes K
$$
et de ce que $\a_{H^D}(-1)\circ \a_H=0$ dans la suite de Hodge-Tate de
$H^D$. Donc, $\zeta_G$ définit un élément de $V_p (G)$.

\paragraph{Cinquième étape : Rigidification }

Il suffit de montrer que $\zeta_G$ est non-nul mais cela est clair. 
On obtient donc d'après le théorème \ref{descDri} un triplet
$(G,\rho_G,\eta_G)\in \M_\infty^\D (K)$. 

\subsection{Les deux applications sont inverses l'une de l'autre}

Cela est moins clair que dans la première description de
l'isomorphisme. 
\\

Partons de $(G,\rho_G,\eta_G)\in \M_\infty^\D (K)$ et soit
$(H,\rho_H,\eta_H)\in \M_\infty (K)$ le triplet associé. Soit 
$(G',\rho_{G'},\eta_{G'})\in\M_\infty^\D (K)$ le triplet associé à
$(H,\rho_H,\eta_H)$. Rappelons que la filtration de $\DG_{\Q}\otimes
K$ 
 définissant
$(G',\rho_{G'})\in \M^\D (K)/\sim$ s'identifie à
$$
\left [ V_p (H)^*\otimes K \twoheadrightarrow \omega_H\otimes K(-1)
\right ] \otimes_L \DH_\Q
$$
via 
$$
V_p (H)^* \otimes \DH_\Q \underset{\text{via }\eta_H}{\xrig{\;\;\sim\;\;}}
\DH_\Q^n \xrig{\;\mathbb{D} (\Delta)\;} \DG_\Q
$$
Rappelons que l'on a un isomorphisme  
$$
\Hom_{\O_D} ( \Gb,\Hb) \iso V_p (H)
$$
et qu'alors l'identification ci-dessus se résume à
$$
\xymatrix@R=2mm{
\DG_\Q \ar[r]^(.24)\sim & \Hom_F \left ( \Hom_{\O_D} (\Gb,\Hb)\unp,
  \DH_\Q\right )
  & \ar[l]_(.34)\sim  V_p (H)^*\otimes \DH_\Q \\
x \ar@{|->}[r] & \left [ f \longmapsto \mathbb{D} (f) (x) \right ]
}
$$
D'après l'exactitude de la suite de Hodge-Tate de $H^D$
$$
\ker \left ( V_p (H)^*\otimes K\twoheadrightarrow \omega_G\otimes
  K(-1) 
\right)\otimes_L \DH_\Q = \{ \;h :V_p (H)\otimes K\ldrt \DH_\Q\otimes K
\; |\; \a_H (x) =0 \impl h(x) =0 \;\}
$$
Mais via l'identification $\Hom_{\O_D} (\Gb,\Hb)\unp \simeq V_p (H)$
l'application $\a_H$ est 
$$
f\longmapsto \left ( \mathbb{D} (f) \otimes Id \right ) ( \zeta_{G*}
(1) )
$$
où $\mathbb{D} (f)\otimes Id : \DG\otimes K\ldrt \DH\otimes K$,
$\zeta_G= \rho_G \circ \eta_G (1)$ et $\zeta_{G*}$ est l'application
induite sur les cristaux. 
On en déduit que la filtration définissant $(G',\rho_{G'})$ sur
$\DG_\Q\otimes K$ est 
$$
\{\, x\in \DG_\Q\otimes K \; |\; \forall f \in \Hom_{D,\ph} \left (
  \DG_\Q,\DH_\Q \right ) \; (f\otimes Id) (\zeta_{G*} (1))=0 \impl
(f\otimes Id) (x) =0 \;\} =\left ( Im (\a_G)^\perp \right )^\perp
$$
qui est donc égal à $Im \a_G $ qui par surjectivité de l'application de
Hodge-Tate de $G$ est la filtration définissant $G$.
Donc $(G',\rho_{G'}) = (G,\rho_G)$.

Il est maintenant aisé de vérifier que les rigidifications des modules
de Tate de $G'$ et $G$ coïncident puisqu'il suffit de vérifier
qu'elles coïncident modulo $p$, dans $\Hom (F/\O_F, \Gb\otimes
\O_K/p\O_K)$ ce qui est immédiat. Donc $(G',\rho_{G'},\eta_{G'}) =
(G,\rho_G, \eta_G)$.
\\

On vérifie de la même façon que l'application composée 
$$
\M_\infty^{\LT} (K) \ldrt \M_\infty^\D (K) \ldrt \M_\infty^{\LT} (K)
$$
est l'identité.

\appendix

\section{Théorèmes de comparaison entiers relatifs pour les groupes
  $p$-divisibles d'après Faltings}

Dans cet appendice on explique les résultats auxquels on peut parvenir
à partir des méthodes de  \cite{Faltings1} pour les périodes
cristallines et de Hodge-Tate des
groupes $p$-divisibles sur les anneaux d'entiers de corps
non-archimédiens. Les démonstrations seront données dans
\cite{Periodes}.

\subsection{Groupes p-divisibles sur les anneaux d'entiers de corps non-archimédiens}

Soit $K|\Qp$ un corps valué complet pour une valuation de rang $1$,
c'est à dire à valeurs dans $\R$, étendant celle de $\Qp$. On note
$\O_{K}$ son anneau des entiers et $k$ son corps résiduel. Fixons
$\O_{K_0}\subset \O_{K}$ un anneau de Cohen. On a des extensions
valuées $K|K_0 |\Qp$ . 
 Le choix de l'anneau de Cohen fixe en particulier une section $\e$ du morphisme $\O_{K}/p\O_{K}
\twoheadrightarrow k$. 
  Si $k$ est parfait $\O_{K_0}\simeq W(k)$.

Par définition un groupe $p$-divisible sur $\spf (\O_K)$ est un
système compatible de groupes $p$-divisibles sur les $(\spec
(\O_{K}/p^n \O_{K}))_{n\geq 1}$ (cf. le début de la section 1 de \cite{Cellulaire}). 

\begin{lemm}
Soit $G$ un groupe $p$-divisible sur $\spf (\O_K)$. 
Son équivalents :
\begin{itemize}
\item $G\otimes_{\O_{K}} \O_{K}/p \O_{K}$ est isogéne à un groupe constant $H$ sur $k$ : $$G\otimes_{\O_{K}} \O_{K}/p \O_{K}\sim H\otimes_{k,\e} 
\O_{K}/p \O_{K}$$
\item Il existe un nombre réel $\l \geq 1$ et $H'$ un groupe p-divisible  sur $k$ tels que 
si $\mathfrak{m}_{K,\l}=\{ x \in K \;|\; v_p (x)\geq \l \;\}$ alors
$$
G\otimes \O_{K}/ \mathfrak{m}_{K,\l} \simeq H'\otimes_{k,\e}  \O_{K}/ \mathfrak{m}_{K,\l}
$$
\item Soit $R\simeq \O_{K_0}[[x_i]]_{i\in I}$ l'anneau universel des
  déformations du groupe $G\otimes k$ sur des $\O_{K_0}$-algèbres locales complètes et $G^{univ}$ la déformation universelle. 
Il existe un morphisme
$$
x : \spf (\O_{K}) \ldrt \spf ( R)
$$
i.e. $x\in  \spf ( R)^{an} (K)$ où $\,^{an}$ désigne la fibre générique au sens des espaces de Berkovich, tel que $G\simeq x^* G^{univ}$
\end{itemize}
\end{lemm}
\begin{proof}    
Elle ne pose pas de problème.
\end{proof}

On remarquera que si les conditions du lemme précédent sont vérifiées
alors nécessairement $H'\simeq G_k$ et on peut choisir $H=G_k$.

\begin{defi}
Un groupe p-divisible satisfaisant les conditions équivalentes du
lemme précédent sera dit isotrivial mod $p$.
\end{defi}

\begin{rema}
En utilisant des morphismes non-continus $\O_{K_0}[[x_i]]_{i\in I}
\ldrt \O_{\Cp}$ on peut construire des groupes $p$-divisibles
non-isotriviaux mod $p$ sur $\O_{\Cp}$.
\end{rema}

\subsection{Théorèmes de comparaison}

On reprend les notations de la section précédente. On fixe un relèvement de Frobenius $\s:  \O_{K_0}\drt  \O_{K_0}$. Soit $G$ un groupe p-divisible sur $\spf (\O_{K})$. 
On note $M$ l'algèbre de Lie de l'extension vectorielle universelle de $G$. Celle-ci est filtrée :
$$
0 \ldrt \omega_{G^D} \ldrt M \ldrt \omega_G^* \ldrt 0
$$
On note $\text{Fil}\, M = \omega_{G^D}$. On note $\mathcal{E}$ le cristal de Messing (covariant)
 de $G\otimes \O_{K}/ p \O_{K}$ sur le gros site cristallin nilpotent
 de \cite{BBM} $$
NCRIS(\spec (\O_{K}/p \O_{K})/ \spec (\O_{K_0}))
$$ 
où $\spec (\O_{K_0}) \ldrt \spec (\O_{K}/p \O_{K})$ est défini via la
section $\e$, 
et  
$$
M_0 = \mathcal{E}_{\O_{K_0}\twoheadrightarrow k}
$$
le module de Dieudonné ``classique'' de $G\otimes k$ si $k$ est
parfait. Il est muni d'une application $\s$-linéaire $\varphi:M_0\ldrt
M_0$ associée au relèvement de Frobenius $\s$.  

Rappelons (\cite{Periodes}) que l'on dispose d'une $\O_{K_0}$-algèbre $A_{cris} (\O_{K})$ augmentée via $\theta : A_{cris} (\O_{K}) \twoheadrightarrow \mathcal{O}_{\widehat{\overline{K}}}$. On note $\varphi $ le Frobenius cristallin sur $A_{cris}$. 

L'évaluation 
$$
E = \mathcal{E}_{ A_{cris} \twoheadrightarrow \mathcal{O}_{\widehat{\overline{K}}}}
$$
est un $A_{cris}$-module libre muni d'un morphisme $\varphi$-linéaire $\varphi : E \ldrt E$ tel que $$E\otimes_{\theta} \mathcal{O}_{\widehat{\overline{K}}} \simeq M\otimes_{\O_{K}} \otimes \mathcal{O}_{\widehat{\overline{K}}}$$ ce qui permet de filtrer $E$ via 
$\text{Fil}\, E = \theta^{-1} (\text{Fil}\, M \otimes \mathcal{O}_{\widehat{\overline{K}}})$. On a donc un $\varphi$-module filtré $(E,\varphi, \text{Fil}\, E)$. 

\begin{theo}[\cite{Periodes}]\label{crop}
Il y a un isomorphisme naturel $\Gal (\overline{K}|K)$-équivariant 
$$
T_p (G) \iso \left ( \text{Fil}\; E \right )^{\varphi =p}
$$
\end{theo}

\begin{theo}[\cite{Periodes}]\label{bouqibouc}
Il y a deux inclusions naturelles strictement compatibles aux filtrations et à l'action de Galois
$$
t E \subset T_p (G)\otimes_{\Zp} A_{cris} \subset E 
$$
Elles sont compatibles aux Frobenius cristallins lorsque l'on munit $T_p (G) \otimes A_{cris}$ de $p\otimes \ph$.
\end{theo}

\begin{rema}
Dans le théorème précédent la filtration est indexée de la façon
suivante 
\begin{eqnarray*}
\Fil^{\;i} \left ( T_p (G)\otimes A_{cris}\right ) &=& T_p (G)\otimes \Fil^{\; i} A_{cris} \\
\forall i\leq -1 \;\Fil^{\;i}\, E &=& E \\
\Fil^{\; 0}\, E &= & \theta^{-1} (\text{Fil} \,M \otimes \mathcal{O}_{\widehat{\overline{K}}}) \\
\forall i\geq 1\; \Fil^{\; i}\, E &=& \Fil^{\; i} A_{cris}. \Fil^{\; 0} 
\end{eqnarray*}
\end{rema}

Lorsque le groupe $p$-divisible $G$ est isotrivial mod $p$ son isocristal
est engendré par ses section horizontales. On en déduit le théorème suivant.

\begin{theo}\label{comp_modulaire}
Supposons de plus que $G$ est isotrivial mod $p$. Il y a alors des isomorphismes 
$$
M_0 \unp \otimes_{K_0} K \simeq M \unp
$$
et 
$$
E \unp \simeq M_0 \unp \otimes_{K_0} B^+_{cris}
$$
comme $\varphi$-modules où $B^+_{cris} = A_{cris} \unp$. 

Il y a donc des inclusions strictement compatibles aux filtrations 
$$
t M_0 \unp \otimes_{K_0} B^+_{cris} \subset V_p (G) \otimes_{\Qp} B^{+}_{cris} \subset M_0 \unp \otimes_{K_0} B^+_{cris} 
$$
 un isomorphisme 
$$
M_0\unp\otimes_{K_0} B_{cris} \simeq V_p (G) \otimes_{\Qp} B_{cris}
$$
avec $B_{cris} =B^+_{cris} \unp$ et un autre isomorphisme
$$
V_p (G) \iso \Fil \left ( M_0\unp \otimes_{K_0} B^+_{cris}\right )^{\varphi = p }
$$
où sur $M_0\otimes B^+_{cris}$, $\varphi= \varphi \otimes \varphi$ et $\Fil$ est la filtration associée 
à celle de $M\unp$ via l'isomorphisme $M_0\unp \otimes_{K_0} K \simeq M\unp$ et 
$\theta$. 
\end{theo}

\subsection{Le déterminant des périodes divisé par $2i\pi$ est une unité
  $p$-adique}

\subsubsection{\'Enoncé et dévissage au cas C.M.} 

Soit $K$ comme précédemment et supposons de plus que son corps
résiduel est algébriquement clos (ce que l'on peut toujours réaliser
quitte à étendre les scalaires). 
Soit $G$ un groupe $p$-divisible modulaire sur $\spf (\O_{K})$ de
  dimension $d$ et hauteur $h$. 
  L'isomorphisme de $B_{cris}$-modules
$$ 
V_p (G) \otimes_{\Qp} B_{cris} \iso M_0 \unp \otimes_{K_0} B_{cris}
$$
du théorème \ref{comp_modulaire} induit un isomorphisme $\ph$-équivariant
$$ \a: 
\det V_p (G) \otimes_{\Qp} B_{cris} \iso ( \det M_0 )\unp \otimes_{K_0} B_{cris}
$$
Plus précisément, il y a une inclusion de $B^+_{cris}$-modules
compatible aux filtrations et Frobenius 
$$
u : T_p (G) \otimes B^+_{cris} \hookrightarrow M_0\otimes_{\O_{K_0}} B^+_{cris}
$$
de conoyau annulé par $t$, d'où une inclusion
$$
\det u : \det T_p (G) \otimes_{\Zp} B^+_{cris} \hookrightarrow \left (
  \det M_0 \right ) \otimes_{\O_{K_0}} B^+_{cris}
$$
 Via l'application de réduction modulo $\Fil^1 B^+_{cris}$ et
l'identification $M_0 \otimes_{\O_{K_0}} \widehat{\overline{K}} \simeq
M \otimes_{\O_K} \widehat{\overline{K}}$ (où $M=\Lie\, E (G)$)
$$
\left ( M_0 \otimes_{\O_{K_0}} B^+_{cris} \right ) / \Fil^1
B^+_{cris}.\left ( M_0 \otimes_{\O_{K_0}} B^+_{cris} \right )
\simeq M \otimes_{\O_K} \widehat{\overline{K}}
$$
et $\Fil\left ( M_0\otimes B^+_{cris}\right ) $ est l'image réciproque de
  $\Fil \, M \otimes  \widehat{\overline{K}}$ (où $\Fil\, M =
  \omega_{G^D}$).
Il existe donc d'après le lemme de Nakayama une base $(e_1,\dots,e_n)$
du $B^+_{cris}$-module $M_0\otimes B^+_{cris}$ telle que 
$$
\Fil \left ( M_0\otimes_{\O_{K_0}} B^+_{cris}\right ) = B^+_{cris}
e_1\oplus \dots \oplus B^+_{cris} e_{n-d-1} \oplus \Fil^1 B^+_{cris}
e_{n-d} \oplus \dots \Fil^1 B^+_{cris} e_n
$$
L'image de $u$ est incluse dans $\Fil\left (M_0\otimes B^+_{cris}\right
)$. On en déduit que 
$$
\det u : \det T_p (G) \otimes B^+_{cris} \hookrightarrow \det M_0
\otimes_{\O_{K_0}} \Fil^d B^+_{cris}
$$
qui définit donc un élément 
$$
\beta \in \left ( \det T_p (G)^{-1} \otimes_{\Zp} \det M_0
  \otimes_{\O_{K_0}} \Fil^d\, B^+_{cris} \right )^{\ph = Id} 
$$
où rappelons que $\ph$ agit par la multiplication par $p$ sur $T_p
(G)$ (cf. théorème \ref{bouqibouc}).
Mais le corps résiduel de $K$ étant algébriquement clos 
$$
(\det M_0, \det \ph) \simeq (\O_{K_0}, p^{h-d} \s)
$$
Donc, 
$$
\beta \in \det (T_p (G))^{-1} \otimes_{\Zp} \left ( \det ( M_0)\right
)^{\ph= p^{h-d}}
\otimes_{\Zp} \left (\Fil^d\,  B^+_{cris} \right )^{\ph = p^d}
$$
Mais puisque  (cf. \cite{Periodes}) 
$$
\left (\Fil^d\,  B^+_{cris} \right )^{\ph = p^d} = \Qp.t^d
$$
Donc
$$
\beta \in \left ( \det T_p (G)\right )^{-1}\otimes \left ( \det
  M_0\right )^{\ph= p^{h-d}} \otimes \Qp.t^d
$$
Les structures entières  $\det T_p (G) \simeq \Zp$, $\det M_0 \simeq \O_{K_0}$ induisent 
une $\Zp$-structure sur $\Qp t^d$. L'élément $\beta$ fournit donc un élément de 
$
\Qp^\times/\Zp^\times. t^d
$

\begin{prop} [Faltings] \label{det_periodes_unite}
 L'élément $\beta$ est entier : $\beta \in \Zp^\times.t^d$. 
\end{prop}
\dem
 Le groupe $G$ étant modulaire on peut le mettre en famille : pour un
 morphisme 
$x:\spf (\O_K)\drt \spf (\O_{K_0} [[x_i]]_{1\leq i\leq d(h-d)})$, $G=x^*G^{univ}$
où $G^{univ}$ désigne la déformation universelle de la fibre spéciale
de $G$.  
 Il existe de plus
un rationnel $\a\in\Q_{>0}$ tel que $\forall i \; v_p ( x^* (x_i))\geq \a$. Le morphisme 
$x$ se factorise donc en 
$$
x: \spf (\O_{K}) \ldrt \mathcal{C} \ldrt \spf (\O_{K_0} [[x_i]]_i)
$$
où $\mathcal{C}=\spf (R)$ est un modèle formel $p$-adique normal sans
$p$-torsion sur $\spf (O_{K_0})$ de la boule rigide formée
des éléments de valuation supérieure ou égale à $\a$. Il y a donc un morphisme de spécialisation
$A_{cris} (R)\twoheadrightarrow A_{cris} (\O_{K})$. 
Sur $\spec (R/pR)$ le groupes $H\otimes R/pR$ est isogéne au groupe définissant l'espace des déformations.
 Donc, 
si $H$ désigne le
groupe $p$-divisible restriction de $G^{univ}$ à  $\mathcal{C}$ il y a un isomorphisme de comparaison 
$$ 
V_p (H) \otimes_{\Qp} B_{cris} (R) \iso M_0 \unp \otimes_{K_0} B_{cris} (R)
$$ 
(cf. \cite{Periodes}). 
On obtient comme précédemment un élément  
$$
\beta'\in  \left ( ( \det T_p (H))^{-1} \otimes_{\Zp} \det M_0 \otimes_{\O_{K_0}} \Fil^d\, B_{cris}^+(R)\right )^{\ph = Id} \simeq \Qp t^d
$$
où la dernière égalité résulte de \cite{Periodes}.
 Et bien sûr, via $A_{cris} (R) \twoheadrightarrow A_{cris} (\O_{K})$, $\beta'\mapsto \beta$. 

L'énoncé du théorème est donc invariant par transport parallèle : on peut transporter $\beta$ parallèlement en n'importe quel point du disque unité $\spf (\O_{K_0}[[x_i]]_i)^{an}$ (quitte à agrandir 
le rayon de la boule).
Le mieux est de choisir un point pour lequel la matrice des périodes est diagonalisable. 
 On peut par exemple prendre un point C.M. ayant multiplication complexe
par $\Z_{p^h}$ où $h$ est la hauteur de $G$. En effet, l'action de $\Z_{p^h}$ permet de diagonaliser la matrice des périodes si l'on prend des bases de vecteurs propres pour cette action. 
 Le résultat est  démontré dans la section suivante
dans  ce cas particulier.
\qed

\subsubsection{Étude des périodes entières des groupes p-divisibles ayant multiplication complexe par un ordre maximal  non-ramifié}\label{CM}

Nous étudions ici les périodes des groupes p-divisible C.M. les plus simples, ceux ayant multiplication
complexe par l'anneau des entiers d'une extension non-ramifiée de $\Qp$. 
Soit donc $G$ un groupe p-divisible de hauteur $h$ et dimension $d$ 
sur $\Z_{p^h}$ muni d'une action  $\iota:\Z_{p^h}\drt \End (G)$. On pourra par exemple prendre lorsque $d=1$ le groupe formel d'exponentielle 
$$
f(T)=\sum_{n\geq 0} \frac{T^{p^{nh}}}{p^n}
$$
qui est bien muni d'une action de $\Z_{p^h}$ puisque $\forall \zeta \in \mu_{p^h-1}\;\; f(\zeta T) =\zeta f(T)$ ou n'importe quel groupe de Lubin-Tate de hauteur $1$ pour l'extension $\Q_{p^h}$.

Notons 
$$
\chi: \Gal(\Qpb |{\Q_{p^h}}) \ldrt \Z_{p^h}^\times
$$
le caractère de Lubin-Tate. Notons $\s \in \Gal (\Q_{p^h}|\Qp)$ le Frobenius et $\forall i\; \chi^{\s^i}= \s^i \circ \chi$. Il résulte de la classification des représentations cristallines abeliennes (Fontaine, cf. \cite{Wintemberger1}) qu'il
existe des entiers $a_i$, $0\leq a_i \leq d$ tels que 
$$
V_p (G)= \prod_{i=0}^{h-1} \chi^{a_i \s^i}
$$
où $\sum_i a_i = d$ et $a_i=\text{rg}_{\Z_{p^h}} \text{Lie} (G)_i$
avec
$$
\Lie G = \bigoplus_{i\in \Z/h\Z} \Lie (G)_i
$$
$\iota (\Z_{p^h})$ agissant sur 
$ \text{Lie} (G)_i$ via $\s^i$. 

Considérons l'application de comparaison 
$$
T_p (G) \otimes_{\Zp} A_{cris} (\Z_{p^h}) \hookrightarrow M_0\otimes_{\Z_{p^h}} A_{cris} (\Z_{p^h})
$$
où $M_0$ est le module de Dieudonné de la fibre spéciale de $G$. 
L'action $\iota$ permet de décomposer 
\begin{eqnarray*}
T_p (G) \otimes_{\Zp} A_{cris} (\Z_{p^h}) &=& \bigoplus_{i\in \Z/ h\Z} \left (T_p (G) \otimes_{\Zp} A_{cris} (\Z_{p^h})\right )_i \\
M_0 &=& \bigoplus_{i\in \Z/h \Z} M_{0,i}
\end{eqnarray*}
où $\iota (\Z_{p^h})$ agit sur la composante indexée par $i$ via $\s^i$. La naturalité du morphisme de comparaison 
implique que celui-ci est somme directe de morphismes 
$$
\left (
T_p (G) \otimes_{\Zp} A_{cris} (\Z_{p^h}) \right )_i \hookrightarrow M_{0,i} \otimes_{\Z_{p^h}} A_{cris} (\Z_{p^h})
$$
Après choix d'une base de $T_p (G)$ et $M_{0,i}$ comme $\Z_{p^h}$-modules ce morphisme est donnée par la 
multiplication par un élément $x_i \in A_{cris}$ bien défini modulo $\Z_{p^h}^\times$. 
D'après la compatibilité stricte aux filtrations du morphisme de comparaison si $a_i=0, x_i\notin \Fil^1 A_{cris}$ et $x_i \in \Fil^1 A_{cris}$ sinon. 
Notons alors $y_i\in \O_{\C_p}$ l'image de $x_i$ dans $\text{Gr}^0 A_{cris}$, resp. $\text{Gr}^1 A_{cris}$. 

\begin{prop}
Supposons $d=1$ et soit $i_0$ l'unique indice tel que $a_{i_0}\neq 0$. Alors,
\begin{eqnarray*}
\forall i<i_0 \; v(y_i) &=& \frac{p^{h+i- i_0}}{p^h-1} \\
\forall i\geq i_0 \; v(y_i) &=& \frac{p^{i-i_0}}{p^h -1}
\end{eqnarray*}
\end{prop}
\dem
La compatibilité  du morphisme de comparaison à l'action de $\Gal(\Qpb
| {\Q_{p^h}})$ est équivalente à ce que 
$$
\forall \tau \in G_{\Q_{p^h}}\;\;\; x_i^\tau =\chi_G (\tau)^{\s^i} x_i 
$$
où $\chi_G$ désigne le caractère galoisien associé à $V_p (G)$. 
Puisque $d=1$, $\chi_G=\chi^{\s^{i_0}}$. Considérons maintenant le lemme suivant : 
\begin{lemm}
Soit $z\in \O_{\C_p}$ tel que $\forall \tau \in \Gal (\Qpb|\Q_{p^h}) \;\; z^\tau = \chi (\tau)^{\s^i} z$. Alors, 
$$
\exists j\in \N\;\; v(z)\in \frac{p^i}{p^{jh} (p^h -1)} + \Z
$$
\end{lemm}
\dem 
Pour $\l \in \Q$ notons $\O_{\C_p,\geq \l}$, resp. $\O_{\C_p,> \l}$, les éléments de valuations supérieure à $\l$, resp. strictement supérieure à $\l$. Notons $q=p^h$. Décrivons l'action de $\Gal(\Qpb|\Q_{p^h})$ sur 
$\O_{\C_p,\geq \l}/\O_{\C_p,> \l} \simeq \overline{\Fp}$. Fixons
$p^\l\in \Qpb$ un élément tel  que si 
$\l=\frac{r}{s}$ avec $r\wedge s=1$ 
 on ait $
(p^\l)^s = p^r$. Soit $\tau_0\in \Gal (\Qpb | \Q_{p^h}(p^\l))$ un
relèvement du Frobenius $x\mapsto x^{p^h}$.
 L'inertie $I_{\Q_{p^h}}$ agit sur $\O_{\C_p,\geq \l}/\O_{\C_p,> \l}$ via le caractère modéré 
$$
\tau \mapsto \frac{\tau (p^\l)}{p^\l} \text{ mod } \O_{\C_p,> \l} \in \overline{\Fp}^\times
$$
Écrivons $z= u p^\l$ où $u\in \O_{\C_p}^\times$. Alors, modulo $\O_{\C_p,> \l}$
\begin{eqnarray*}
\forall \tau \in I_{\Q_{p^h}}\; \forall k\in \Z\; \; \frac{\tau_0^k\tau (z)}{z} &\equiv &\frac{\tau_0^k (u)}{u} . \frac{\tau_0^k \tau (p^\l)}{p^\l}  \;\text{ car } \frac{\tau (u)}{u} \equiv 1  \\
&\equiv & u^{q^k-1} \left ( \frac{\tau (p^\l)}{p^\l} \right )^{q^k} \;\text{ car } \tau_0 (p^\l)=p^\l
\end{eqnarray*}
Quant au caractère $\chi$ il vérifie la congruence 
$$
\forall \tau'\in \Gal (\Qpb| \Q_{p^h})\;\; \chi ( \tau') \equiv \frac{\tau' (p^{\frac{1}{q-1}})}{p^{\frac{1}{q-1}}} \in \mu_{q-1}
$$
L'hypothèse du lemme implique que
$$
\forall k\in \Z\; \forall \tau\in I_{\Q_{p^h}} \;\; u^{q^k-1} \frac{\tau ( p^{q^k \l-\frac{p^i}{q-1}})}{
 p^{q^k \l-\frac{p^i}{q-1}}} \equiv 1 
$$
ce qui implique d'abord avec $\tau= Id$ que $\bar{u}\in {\mathbb{F}_q}^\times$ et que l'on peut donc supposer que $u=1$ dans
la congruence ci-dessus. De plus, pour $\mu \in\Q$, le caractère modéré
\begin{eqnarray*}
I_{\Q_{p^h}}&\ldrt & \overline{\Fp}^\times \\
\tau & \longmapsto & \frac{\tau(p^\mu)}{p^\mu}
\end{eqnarray*}
est trivial ssi $\mu\in \Z [ \frac{1}{p}]$. Donc, $\forall k\in \Z\; q^k \l -\frac{p^i}{q-1} \in \Z[\frac{1}{q}]$ ce qui implique facilement le lemme.
\qed

Il résulte du lemme précédent que 
\begin{eqnarray*}
\forall i>i_0\; \exists j\in \N\; v_p (y_i)\in \frac{p^{i-i_0}}{q^j ( q-1)} + \Z \\
\text{ et } \forall i<i_0\; \exists j\in \N\; v_p (y_i)\in \frac{p^{h+i - i_0}}{q^j (q-1)}+\Z 
\end{eqnarray*}
Quant à $y_{i_0}$, étant donné que dans $\text{Gr}^1 A_{cris}$ $t$ se transforme via $N_{\Q_{p^h}/\Qp} \circ \chi$ et que $v(t)=\frac{1}{p-1}$, 
$$
\exists j\in \N\; v_p (y_{i_0})\in \frac{1}{q^j (q-1)} + \Z 
$$
Le morphisme de comparaison possède un quasi-inverse tel qu'avec composition avec celui-ci on obtienne 
la multiplication par $t$ qui est de valuation $1/(p-1)$. On en déduit que 
\begin{eqnarray*}
\forall i \leq i_0 \;\; 0 \leq v(y_i) \leq \frac{p^{i-i_0}}{q-1}  \\
\forall i <i_0 \; \;  0 \leq v(y_i) \leq \frac{p^{h+i-i_0}}{q-1}
\end{eqnarray*}
Remarquons maintenant que 
$$
\prod_i x_i = \beta t 
$$
où $\beta \in \Zp$, et 
$$
v(\beta) + \frac{1}{p-1} = \sum_i v_p (y_i) \leq \frac{1}{p-1}
$$
et donc $v(\beta)=0$ i.e. $\beta \in \Zp^\times$ et les valuations des $y_i$ sont celle annoncées. 
\qed

Il résulte donc de la démonstration de la proposition précédente que 

\begin{coro}
L'énoncé du théorème \ref{det_periodes_unite} est vrai dans le cas de dimension $1$ et de multiplication complexe par $\Z_{p^h}$· 
\end{coro}

Attaquons maintenant le cas de dimension $d$ quelconque. Commençons par remarquer qu'il résulte de l'étude du cas de dimension $1$  et du théorème de Tate, $H^0 (G_{\Q_{p^h}},\C_p)=\Q_{p^h}$, 
 que 

\begin{coro}
Pour tout entier $i$ compris entre $0$ et $h-1$, 
$$
H^0 (G_{\Q_{p^h}}, \C_p (\chi^{-\s^i}))=\Q_{p^h}.z_i
$$
où $z_i\in \C_p$ est un élément de valuation $\dpt{\frac{p^i}{p^h-1}}$. 
\end{coro}

Il en résulte que l'on a une majoration des valuations des éléments $y_i$ (les périodes partielles).
Sachant que $\prod_i y_i = \beta $ on en déduit le résultat facilement. 
\qed

On renvoie également à l'article \cite{Colmez1}.

\bibliographystyle{plain}
\bibliography{biblio}
\end{document}